\documentclass[12pt]{amsart}

\usepackage{amssymb}
\usepackage{amsthm}
\usepackage{amscd}

\usepackage{pifont}

\usepackage{amsmath}
\usepackage{latexsym}

\usepackage{graphicx}
\usepackage{amsfonts}

\usepackage[all]{xy}

\theoremstyle{plain}
\newtheorem{lemma}{Lemma}[subsection]
\newtheorem{prop}[lemma]{Proposition}
\newtheorem{thm}[lemma]{Theorem}
\newtheorem{cor}[lemma]{Corollary}
\newtheorem{aplemma}{Lemma~A.\hspace{-1.5mm}}
\newtheorem{approp}{Proposition~A.\hspace{-1.5mm}}
\newtheorem{apthm}{Theorem~A.\hspace{-1.5mm}}
\newtheorem{apcor}{Corollary~A.\hspace{-1.5mm}}
\newtheorem{intthm}{Theorem}

\newtheorem{conj}[lemma]{Conjecture}

\theoremstyle{definition}

\newtheorem{rem}{Remark}
\newtheorem{rema}[lemma]{Remark}

\newtheorem{remb}{Remark}

\newtheorem{defi}[lemma]{Definition}
\newtheorem{exa}[lemma]{Example}
\newtheorem{aprem}{Remark~A.\hspace{-1.5mm}}
\newtheorem{apdefi}{Definition~A.\hspace{-1.5mm}}
\newcommand{\bde}{\begin{defi}}
\newcommand{\ede}{\end{defi}\vspace{1mm}}
\newcommand{\ble}{\begin{lemma}}
\newcommand{\ele}{\end{lemma}}
\newcommand{\bpr}{\begin{prop}}
\newcommand{\epr}{\end{prop}}
\newcommand{\bt}{\begin{thm}}
\newcommand{\et}{\end{thm}}
\newcommand{\bco}{\begin{cor}}
\newcommand{\eco}{\end{cor}}
\newcommand{\bre}{\begin{rem}}
\newcommand{\ere}{\end{rem}}
\newcommand{\brea}{\begin{rema}}
\newcommand{\erea}{\end{rema}\vspace{1mm}}
\newcommand{\breb}{\begin{remb}}
\newcommand{\ereb}{\end{remb}\vspace{1mm}}
\newcommand{\bex}{\begin{exa}}
\newcommand{\eex}{\end{exa}}
\newcommand{\bpf}{\begin{proof}}
\newcommand{\epf}{\end{proof}\vspace{1mm}}

\newcommand{\bade}{\begin{apdefi}}
\newcommand{\eade}{\end{apdefi}}
\newcommand{\bale}{\begin{aplemma}}
\newcommand{\eale}{\end{aplemma}}
\newcommand{\bapr}{\begin{approp}}
\newcommand{\eapr}{\end{approp}}
\newcommand{\bat}{\begin{apthm}}
\newcommand{\eat}{\end{apthm}}
\newcommand{\baco}{\begin{apcor}}
\newcommand{\eaco}{\end{apcor}}
\newcommand{\bare}{\begin{aprem}}
\newcommand{\eare}{\end{aprem}}

\newcommand{\be}{\begin{enumerate}}
\newcommand{\ee}{\end{enumerate}}
\newcommand{\bcd}{\[\begin{CD}}
\newcommand{\ecd}{\end{CD}\]}
\newcommand{\bit}{\begin{itemize}}
\newcommand{\eit}{\end{itemize}}
\newcommand{\bq}{\begin{quote}}
\newcommand{\eq}{\end{quote}}
\newcommand{\ba}{\begin{array}}
\newcommand{\ea}{\end{array}}

\newcommand{\mcD}{\mathcal{D}}
\newcommand{\mcE}{\mathcal{E}}
\newcommand{\mcF}{\mathcal{F}}
\newcommand{\mcG}{\mathcal{G}}
\newcommand{\mcH}{\mathcal{H}}

\newcommand{\mcK}{\mathcal{K}}
\newcommand{\mcL}{\mathcal{L}}

\newcommand{\mcN}{\mathcal{N}}
\newcommand{\mcO}{\mathcal{O}}

\newcommand{\mcT}{\mathcal{T}}

\newcommand{\mcV}{\mathcal{V}}

\newcommand{\mbC}{\mathbb{C}}

\newcommand{\mbH}{\mathbb{H}}

\newcommand{\mbN}{\mathbb{N}}

\newcommand{\mbP}{\mathbb{P}}

\newcommand{\mbR}{\mathbb{R}}

\newcommand{\mbX}{\mathbb{X}}
\newcommand{\mbY}{\mathbb{Y}}
\newcommand{\mbZ}{\mathbb{Z}}

\newcommand{\mfC}{\mathfrak{C}}
\newcommand{\mfD}{\mathfrak{D}}

\newcommand{\mfM}{\mathfrak{M}}
\newcommand{\mfN}{\mathfrak{N}}
\newcommand{\mfO}{\mathfrak{O}}

\newcommand{\mfR}{\mathfrak{R}}
\newcommand{\mfS}{\mathfrak{S}}

\newcommand{\mfb}{\mathfrak{b}}

\newcommand{\mfg}{\mathfrak{g}}
\newcommand{\mfh}{\mathfrak{h}}

\newcommand{\mfl}{\mathfrak{l}}

\newcommand{\mfo}{\mathfrak{o}}
\newcommand{\mfp}{\mathfrak{p}}

\newcommand{\mfs}{\mathfrak{s}}

\newcommand{\msC}{\mathscr{C}}
\newcommand{\msD}{\mathscr{D}}

\newcommand{\msN}{\mathscr{N}}

\newcommand{\msP}{\mathscr{P}}

\newcommand{\msX}{\mathscr{X}}
\newcommand{\msY}{\mathscr{Y}}

\newcommand{\migi}{\rightarrow}
\newcommand{\longmigi}{\longrightarrow}

\newcommand{\isom}{\stackrel{\sim}{\migi}}

\newcommand{\migiincl}{\hookrightarrow}

\newcommand{\migisurj}{\twoheadrightarrow}

\newcommand{\mr}{\mathrm}
\newcommand{\hidden}[1]{\,}

\newcommand{\SSP}{\vspace{0mm}}
\newcommand{\LSP}{\vspace{0mm}}

\newcommand{\LL}{d}
\newcommand{\LLL}{d}

\usepackage {cancel}
\usepackage{mathrsfs}

\begin{document}

\title[Cyclic coverings and ordinariness of dormant opers]{Cyclic \'{e}tale coverings of generic curves\\ and ordinariness of dormant opers}
\author{Yasuhiro Wakabayashi}
\date{}
\markboth{Yasuhiro Wakabayashi}{}
\maketitle
\footnotetext{Y. Wakabayashi: Department of Mathematics, Tokyo Institute of Technology, 2-12-1 Ookayama, Meguro-ku, Tokyo 152-8551, JAPAN;}
\footnotetext{e-mail: {\tt wkbysh@math.titech.ac.jp};}
\begin{abstract}
The ordinariness of elliptic curves is essential in proving various expected properties of elliptic curves in positive characteristic and can be extended to algebraic curves of arbitrary genus. The present paper deals with another kind of extension, i.e., ordinariness of dormant $\mathfrak{sl}_2$-opers, or more generally, dormant $\mathfrak{sl}_n$-opers. We prove that, for dormant $\mathfrak{sl}_2$-opers on elliptic curves, this notion is essentially equivalent to the classical ordinariness. Moreover,  the main result of the present paper asserts that the pull-back of an ordinary dormant $\mathfrak{sl}_n$-oper on a general curve by a cyclic covering is ordinary whenever the order of its Galois group is prime to the characteristic of the base field. This result may be regarded as an analogue of a result by S. Nakajima for ordinary algebraic curves. 
 \end{abstract}
\tableofcontents

\section{Introduction} \label{wS31}

\LSP
\subsection{} \label{S01}

A {\it (complex) projective structure} on a Riemann surface is a maximal system of local coordinates  in the analytic topology 
modeled on the Riemann sphere $\mbC \sqcup \{ \infty \}$
such that on any two overlapping coordinate patches, the change of coordinates is described as a M\"{o}bius transformation.
An important consequence of the uniformization theorem  is that any closed hyperbolic  Riemann surface $X_\mbC$ is isomorphic  to a quotient of the $1$-dimensional  complex hyperbolic space $\mbH^1_\mbC$ by a torsion-free discrete subgroup of $\mr{SU}(1,1) \cong \mr{SL}_2 (\mbR)$.
This implies that by collecting various local inverses of a universal covering map,
we have a  projective structure on $X_\mbC$.
This projective structure  is canonical because
the Poincar\'{e} metric on $\mbH^1_\mbC$ is invariant under the action by  
 its monodromy group,   so  it  brings  the hyperbolic metric on $X_\mbC$. 
In this way,  such an  additional structure provides  us a foundation for discussing  an enriched  geometry extending  the underlying Riemann surface.

To develop the geometry in characteristic  $p$ based on the theory of complex  projective structures, 
we deal with 
  algebraic curves  equipped with a dormant $\mfs \mfl_2$-oper, which are 
called   {\bf dormant curves} (cf. Definition \ref{wD3}, (i)).
Here, 
let $k$ be  an algebraically closed field of characteristic $p>2$ and $X$  a connected  smooth proper curve over $k$.
A {\bf dormant $\mfs \mfl_2$-oper}  on $X$ is 
defined as a certain $p$-flat $\mr{PGL}_2$-bundle on $X$ equipped with a Borel reduction satisfying a strict form of Griffiths  transversality (cf. Definitions \ref{wD1}, (i)).
We refer the reader to ~\cite[Chap.\,2, \S\,2.1, Definition 2.1, (i), and Chap.\,3, \S\,3.4, Definition 3.15]{Wak5} for the definition of a dormant $\mfg$-oper for  a general semisimple Lie algebra $\mfg$.
Dormant $\mfs \mfl_2$-opers 
 have provided a rich story
under the identifications with equivalent realizations, e.g., Frobenius-projective structures and projective connections with a full set of solutions (cf. ~\cite{Hos2}, ~\cite{Mzk2}, ~\cite{Wak}, ~\cite{Wak2}, ~\cite{Wak3}, ~\cite{Wak10}).
By considering the correspondence with  Frobenius-projective structures (which are certain maximal systems of \'{e}tale coordinate charts),
these objects can be regarded as  analogues of complex projective structures on Riemann surfaces.
In particular, some of the results and observations concerning dormant $\mfs \mfl_2$-opers 
 have been obtained  in analogy with complex projective structures.
For example,   a recent result by the author (cf.   ~\cite{Wak10}) shows that 
a dormant $\mfs \mfl_2$-oper   produces  a Frobenius-constant quantization on a certain space in exactly the same way as a  complex projective structure produces a quantization on the corresponding space over the underlying Riemann surface.
In this way, dormant curves are treated as fundamental objects in our study,
based on their similarity to Riemann surfaces with a projective structure, such as  the canonical one.

\LSP
\subsection{} \label{S01}

In  the present paper,  we focus on a certain class of  dormant curves, which we call {\it ordinary dormant curves}  (cf. Definition \ref{wD2}, (i)).
(The notion of ordinariness was  extended in 
 ~\cite[\S\,6, Definition   6.7.1]{Wak11}  to Frobenius-Ehresmann structures with  an arbitrary local model.) 
Here, recall that an algebraic curve is called ordinary in the classical sense if its $p$-rank is the maximum possible, equal to its  genus.
(See ~\cite{BK}, ~\cite{IR} and ~\cite{Maz} for more general definitions.)
This classical ordinariness is essential in proving  various expected properties of algebraic varieties in positive characteristic.
Indeed, ordinary elliptic curves, or more generally, ordinary Abelian varieties, 
admit canonical Serre-Tate liftings to characteristic $0$,  and 
comparison theorems in 
 $p$-adic Hodge theory can be more easily established for such varieties. 
Also,  the ordinariness condition for algebraic  curves
is related to obstructions to lifting  line bundles trivialized by the relative Frobenius morphisms.

On the other hand, the ordinariness of dormant curves is a kind of nonabelian generalization of the classical definition and 
related to obstructions to lifting  dormant $\mfs \mfl_2$-opers.
In particular, this notion has a direct relationship with the (local and global) study of 
 the moduli space of dormant curves, e.g.,  counting all possible dormant $\mfs \mfl_2$-opers on a given curve.
(The ordinary locus of this moduli space  is interesting in its own right; see ~\cite{Wak10} for a work from a viewpoint  of symplectic geometry.
Also, see Remark \ref{wR29} for the relationship
with 
 the hyperbolically ordinariness of nilpotent indigenous bundle discussed in ~\cite{Mzk2}.)
We say that a smooth curve is {\it dormant-ordinary} if any dormant $\mfs \mfl_2$-oper on this curve is ordinary (cf. Definition \ref{wD2}, (ii)).
For example, 
an elliptic (i.e., genus-$1$) curve is  dormant-ordinary  precisely when  it is ordinary in the classical sense.
Moreover, any dormant $\mfs \mfl_2$-oper on an ordinary elliptic curve comes from a Frobenius-trivialized line bundle.
These imply that the two kinds of  ordinariness are essentially equivalent in this case. 
However, at the time of writing the present paper, 
the author does not know much about the difference between these notions for hyperbolic curves.
Therefore, the question of what properties they induce  in common will naturally arise.

\begin{center}
\begin{picture}(400,120)

\put(15,90){\fbox{$\begin{matrix}
\text{ordinary} \\
\text{elliptic curves} 
\end{matrix}$}}

\put(5,10){\fbox{$\begin{matrix}
\text{dormant-ordinary} \\
\text{elliptic curves} 
\end{matrix}$}}

\put(285,90){\fbox{$\begin{matrix}
\text{ordinary} \\
\text{hyperbolic curves} 
\end{matrix}$}}

\put(285,10){\fbox{$\begin{matrix}
\text{dormant-ordinary} \\
\text{hyperbolic curves} 
\end{matrix}$}}

\put(115,50){Higher Genus Generalizations!}

\put(50,30){\huge{\rotatebox{90}{$\Longleftrightarrow$}}}


\put(110,15){\vector(4,0){160}}
\put(110,92){\vector(4,0){160}}



\end{picture}
\vspace{0mm}
\end{center}

\LSP
\subsection{} \label{S01}

To consider  this question, 
 we  recall  here a result by S. Nakajima concerning  the relationship between
the classical ordinariness and coverings of an algebraic curve (cf. ~\cite[\S\,4, Theorem 2]{Nak}).
Let $Y \migi X$ be a Galois covering between connected proper smooth (hyperbolic) curves with Galois group $G$.
Suppose that  $X$ is  general  in the moduli space of curves
and  $G$ is
a cyclic group of degree prime to $p$.
Then, Nakajima's theorem asserts 
that $Y$ is ordinary.
See ~\cite[\S\,1,  Theorem]{Bou},  ~\cite[Theorem 1.1]{OP}, and ~\cite[Theorem 14]{Ray1}  for more general results.
The aim of this present paper is to prove an analogous result of this theorem for dormant curves, described as follows.

\SSP 
\begin{intthm}  \label{EW1}
Let $\mbX$ and $\mbY$ be dormant hyperbolic curves and $w : \mbY \migi \mbX$ an \'{e}tale covering (cf. Definition \ref{wD3}, (ii)).
\begin{itemize}
\item[(i)]
If $\mbY$ is ordinary, then $\mbX$ is ordinary.
\item[(ii)]
If $\mbX$ is ordinary and general, and $w$ is a cyclic \'{e}tale covering with $p \nmid \sharp (\mr{Gal}(\mbY/\mbX))$, then $\mbY$ is ordinary.
(Here,  we recall that the moduli space classifying  dormant curves of a fixed genus  is irreducible as proved in ~\cite[Chap.\,II, \S\,2.3, Theorem 2.8]{Mzk2}; thus, it makes sense to speak of a ``{\it general}" dormant curve, i.e.,  a dormant curve that determines  a point of this space  lying  outside some fixed closed substack.)
\end{itemize}
\end{intthm}
\SSP

Moreover, we prove the following  generalized assertion for  dormant $\mfs \mfl_n$-opers.
(The ordinariness of dormant $\mfs \mfl_n$-opers on a pointed curve will be defined in Definition \ref{wD043}.)

\SSP
\begin{intthm} 
  \label{y019}
  Let us fix an integer $n$ satisfying either  ``$n=2$ and $2<p$" or ``$2< 2n <p$".
Also, let $\msX$ be a  pointed smooth  hyperbolic  curve over $k$,  $(\msY, w)$ an \'{e}tale covering  (cf. \S\,\ref{S12}) of $\msX$, and $\mcE^\spadesuit$ a dormant $\mfs \mfl_n$-oper on $\msX$.
\begin{itemize}
\item[(i)]
Suppose that  the pull-back $w^*(\mcE^\spadesuit)$ (cf. \S\,\ref{wSdS2}) of $\mcE^\spadesuit$ is ordinary.
Then, $\mcE^\spadesuit$ is  ordinary.
\item[(ii)]
Suppose that $\msX$ is general, $(\msY, w)$ is a cyclic  \'{e}tale covering  with $p \nmid \sharp (\mr{Gal}(\msY/\msX))$, and $\mcE^\spadesuit$ is ordinary.
Then, the pull-back $w^*(\mcE^\spadesuit)$ is ordinary.
\end{itemize}
\end{intthm}
\SSP

To prove the above theorem, we substantially apply various results and discussions  of     ~\cite{Wak5},  in which the author developed the theory 
  of (dormant) opers on pointed stable curves from the viewpoint of logarithmic geometry.
A key ingredient in the proof of 
assertion (i)
is 
the cohomological  description of the deformation space of a dormant $\mfs \mfl_n$-oper
studied in ~\cite[Chap.\,6, \S\,6.5.2]{Wak5}.
If the underlying curve of $\msX$ is smooth,  then  
we can also obtain  this assertion as an application of  
  the correspondence between the moduli space of dormant $\mfs \mfl_n$-opers and a certain Quot scheme established in 
~\cite[Chap.\,9, \S\,9.2, Proposition 9.4]{Wak5} (or ~\cite[\S\,5.4, Proposition  5.4.2]{JP},  
~\cite[\S\,4, Proposition 4.3]{Wak}); indeed,  by this correspondence, the problem follows from  
well-known generalities on the deformation theory of Quot schemes.
Moreover, 
in the discussion on assertion (ii),
we use 
  deformation and degeneration techniques  of the underlying curves
  to 
   reduce the problem  to the case of low genus. 
Thus, the proof is completed by a slightly improved argument regarding dormant $\mfs \mfl_n$-opers on a pointed projective line in ~\cite[\S\,8]{Wak5}.

Finally, we  remark that 
the ordinariness of dormant curves (or more generally, dormant $\mfs \mfl_n$-opers) is essential in the proofs of Joshi's conjecture described in ~\cite{Wak} and ~\cite{Wak5}.
Indeed, we examined  this property  to lift relevant moduli spaces to characteristic $0$ and then applied  a well-known formula for computing the Gromov-Witten invariant of a certain type of Quot scheme over $\mbC$; this formula enables us to 
compute  the number of dormant $\mfs \mfl_n$-opers on a general curve.
For example,  if a hyperbolic curve $X$ over $k$ is both ordinary and  dormant-ordinary, then the main result of ~\cite{Wak} gives an explicit number of dormant $\mfs \mfl_2$-opers on $X$.

\LSP
\subsection{Notation and Conventions} \label{S9025}

 Throughout the present paper,  we fix an odd prime $p$ and 
an algebraically closed field $k$ of characteristic $p$.
All schemes appearing in the present paper are assumed to be locally noetherian.
 
 For each integer $n>1$,
 we denote by  $\mr{PGL}_n$ the rank $n$ projective linear  group over $k$ and by 
 $B_n$ the Borel subgroup of $\mr{PGL}_n$ consisting of elements represented by invertible upper triangular matrices. 
 Write $\mfb_n$ for the Lie algebra of $B_n$.
If $p \nmid n$, then
 the Lie algebra of $\mr{PGL}_n$ may be naturally  identified with $\mfs \mfl_n$, i.e., the Lie algebra  consisting of $n \times n$ matrices over $k$ with vanishing trace.

 Let $f : X \migi S$ be a morphism of $k$-schemes.
Denote by $f^{(1)} : X^{(1)} \migi S$
the  {\it Frobenius twist of $X$ over $S$}, i.e.,  the base-change
 of $f$ via
 the absolute Frobenius morphism $F_S : S \migi S$ of $S$.
Also, 
denote by $F_{X/S} : X \migi X^{(1)}$ the {\it relative Frobenius morphism of $X$ over $S$}.

A {\it log scheme} means, in the present paper,  a scheme equipped with a logarithmic structure in the sense of Fontaine-Illusie (cf.  ~\cite{KaKa}, (1.2)).
For a log scheme (resp., a morphism of log schemes) indicated, say, by the notation ``$X^\mr{log}$" (resp., ``$f^\mr{log}$"), we shall use the notation ``$X$" (resp., ``$f$")  
to denote
 the underlying scheme (resp., the underlying morphism of schemes).

Given a morphism of $k$-log schemes $f^\mr{log} : X^\mr{log} \migi S^\mr{log}$,
we shall write $\mcT_{X^\mr{log}/S^\mr{log}}$ for the sheaf of logarithmic derivations on $X^\mr{log}$ relative to $S^\mr{log}$ and $\Omega_{X^\mr{log}/S^\mr{log}}$ for  its dual, i.e., the sheaf of logarithmic $1$-form on $X^\mr{log}$ relative to $S^\mr{log}$ (cf. ~\cite[(1.7)]{KaKa}).
Next, let $G$ be a  connected smooth affine  algebraic group over $k$ with Lie algebra $\mfg$ and $\pi : \mcE \migi X$ a  $G$-bundle on $X$.
Given a $k$-vector space $\mfh$ equipped with a $G$-action, we write
$\mfh_G$ for the vector bundle on $X$ associated to the affine space $\mcE \times^G \mfh \left(= (\mcE \times_k \mfh)/G \right)$ over $X$.
In particular, by considering 
the case of the $k$-vector space $\mfg$ equipped with the $G$-action given by the adjoint representation $G \migi \mr{GL}(\mfg)$, we obtain the {\it adjoint bundle} $\mfg_\mcE$ associated with $\mcE$.

Let $G$ and $\mcE$ be as above.
By pulling-back the log structure on $X^\mr{log}$ via $\pi$,
one may obtain a log structure on $\mcE$; we denote the resulting log scheme by $\mcE^\mr{log}$.
The $G$-action on $\mcE$ induces  a $G$-action on the direct image $\pi_*(\mcT_{\mcE^\mr{log}/S^\mr{log}})$.
Hence, we obtain the subsheaf of $G$-invariant sections $\widetilde{\mcT}_{\mcE^\mr{log}/S^\mr{log}} := \pi_*(\mcT_{\mcE^\mr{log}/S^\mr{log}})^G \left(\subseteq \pi_*(\mcT_{\mcE^\mr{log}/S^\mr{log}}) \right)$. 
The differential of $\pi$ induces a short exact sequence of $\mcO_X$-modules
\begin{align} \label{wE155}
 0 \longmigi \mfg_\mcE \longmigi \widetilde{\mcT}_{\mcE^\mr{log}/S^\mr{log}} \xrightarrow{d_{\mcE}}\mcT_{X^\mr{log}/S^\mr{log}} \longmigi 0.
 \end{align}
 An {\it $S^\mr{log}$-connection} on $\mcE$ is, by definition,  an $\mcO_X$-linear morphism $\nabla : \mcT_{X^\mr{log}/S^\mr{log}} \migi \widetilde{\mcT}_{\mcE^\mr{log}/S^\mr{log}}$
 satisfying $d_\mcE \circ \nabla = \mr{id}_{\mcT_{X^\mr{log}/S^\mr{log}}}$.
 If $G = \mr{GL}_n$ ($n>0$) and $\mcF$ denotes  the  vector bundle corresponding to $\mcE$, 
 then  this  definition is equivalent to the classical notion of an $S^\mr{log}$-connection on $\mcV$, i.e.,  
 an
$f^{-1}(\mcO_S)$-linear morphism
\begin{equation}
  \nabla : \mcF \migi \Omega_{X^\mr{log}/S^\mr{log}} \otimes \mcF
  \end{equation}
 satisfying  
 $\nabla (a \cdot v) =  d a \otimes v + a \cdot \nabla (v)$
 for any local sections $a \in \mcO_X$, $v \in \mcF$.
 In this situation, we will not distinguish between these notions of connection.
Denote by $\nabla_0$ the trivial $S^\mr{log}$-connection on $\mcO_X$, i.e., the  universal logarithmic derivation $\mcO_X \migi \Omega_{X^\mr{log}/S^\mr{log}}$.

 If $f : X \migi S$ is  a smooth morphism of schemes over $k$ and $\mcG$ is a vector bundle on $X^{(1)}$, then one may define an $S$-connection 
 \begin{align} \label{wE44}
 \nabla_\mcG^\mr{can} : F^*_{X/S}(\mcG) \migi \Omega_{X/S} \otimes F_{X/S}^*(\mcG)
 \end{align}
  on the pull-back $F_{X/S}^*(\mcG)$  determined uniquely by the condition that the sections of the subsheaf $F_{X/S}^{-1}(\mcG)$ are horizontal.
(It is well-known that the assignments $\mcG \mapsto (F^*_{X/S}(\mcG), \nabla_\mcG^\mr{can})$ and $(\mcF, \nabla) \mapsto \mr{Ker}(\nabla)$ define an equivalence of categories between the category of vector bundles on $X^{(1)}$ and the category of vector bundles on $X$ equipped with a flat  $S$-connection with vanishing $p$-curvature. See ~\cite[\S\,5]{Kal}, ~\cite[\S\,1.2]{Og},  or  \cite[Chap.\,3, \S\,3.3, Definition 3.8]{Wak5} for  the definition of $p$-curvature).

Let $\mcE$ be a $G$-bundle on $X$ and $\nabla$ an $S^\mr{log}$-connection on $\mcE$.
Then, we shall write $\nabla^\mr{ad}$ for the $S^\mr{log}$-connection $\mfg_\mcE \migi \Omega_{X^\mr{log}/S^\mr{log}} \otimes \mfg_\mcE$ on the adjoint bundle $\mfg_\mcE$   induced by $\nabla$ via change of structure group by the adjoint representation $G \migi \mr{GL}(\mfg)$.

 \LSP
\subsection*{Acknowledgements} 
We would like to thank Professor Yuichiro Hoshi for  his helpful comments.
Also, we are grateful for the many constructive conversations we had with {\it dormant curves}, who live in the world of mathematics!
Our work was partially supported by Grant-in-Aid for Scientific Research (KAKENHI No. 18K13385, 21K13770).

\vspace{10mm}
\section{Dormant curves and Classification of dormant elliptic curves} \SSP

 In this section, we define the notion of a dormant curve and
 classify  (ordinary) dormant curves of genus $1$.

\LSP
\subsection{Dormant curves} \label{wS11}

 Let $X$ be a connected smooth proper curve over $k$.

\SSP
\bde \label{wD1}
\begin{itemize}
\item[(i)]
Let us consider  a triple  $\mcE^\spadesuit := (\mcE, \nabla, \mcE_B)$ consisting of
a $\mr{PGL}_2$-bundle $\mcE$ on $X$, 
  a $k$-connection $\nabla$  on
  $\mcE$, and a $B_2$-reduction $\mcE_B$ of $\mcE$ (i.e., a $B_2$-bundle $\mcE_B$ together with an isomorphism  of $\mr{PGL}_2$-bundles $\mcE_B \times^{B_2} \mr{PGL}_2 \isom \mcE$).
 We shall say that $\mcE^\spadesuit$ is an {\bf  $\mfs \mfl_2$-oper}  on $X$
     if  
the composite
\begin{align} \label{wE56}
\mr{KS}_{\mcE^\spadesuit} : \mcT_{X/k} \migi \widetilde{\mcT}_{\mcE/k} \migisurj \widetilde{\mcT}_{\mcE/k}/\widetilde{\mcT}_{\mcE_B/k}
\end{align}
is an isomorphism.
\item[(ii)]
Let us consider  a triple   $\mcF^\heartsuit := (\mcF, \nabla, \mcL)$ consisting of 
 a rank $2$ vector bundle on $X$, a $k$-connection $\nabla$  on $\mcF$, and  a line subbundle  $\mcL$ of $\mcF$.
 We shall say that $\mcF^\heartsuit$  is a {\bf $\mr{GL}_2$-oper} on $X$ if   the $\mcO_X$-linear morphism defined as  the composite
 \begin{align} \label{wE55}
\mr{KS}_{\mcF^\heartsuit} :  \mcL \migiincl \mcF \xrightarrow{\nabla} \Omega_{X/k} \otimes \mcF \migisurj \Omega_{X/k} \otimes (\mcF/\mcL)
 \end{align}
 is an isomorphism.
\item[(iii)]
We shall say that an $\mfs \mfl_2$-oper $\mcE^\spadesuit := (\mcE, \nabla, \mcE_B)$ (resp., a $\mr{GL}_2$-oper $\mcF^\heartsuit := (\mcF, \nabla, \mcL)$) is {\bf dormant} if $\nabla$ has vanishing $p$-curvature.

\end{itemize}
\ede
 
 \SSP
 One can define the notion of an isomorphism between two $\mfs \mfl_2$-opers (resp., $\mr{GL}_2$-opers), but we will omit the details of the definition.
 (It is well-known that  an isomorphism
 between given two $\mfs \mfl_2$-opers is, if it exists, 
   uniquely determined. See ~\cite[Chap.\,2, \S\,2.2,  Proposition 2.9]{Wak5}).

Let  $\mcF^\heartsuit := (\mcF, \nabla, \mcL)$ be a  dormant $\mr{GL}_2$-oper  on $X$ and $\msN := (\mcN, \nabla_\mcN)$   a line bundle on $X$ equipped  with a  $k$-connection with  vanishing $p$-curvature.
Denote by $\nabla \otimes \nabla_\mcN$ the $k$-connection on 
the tensor product $\mcF \otimes \mcN$ given by $(\nabla \otimes \nabla_\mcN) (a \otimes b) = \nabla (a) \otimes b + a \otimes \nabla_\mcN (b)$ for any local sections $a \in \mcF$, $b \in \mcN$.
Then, the collection 
\begin{align}
\mcF^\heartsuit_{\otimes \msN} := (\mcF \otimes \mcN, \nabla \otimes \nabla_\mcN, \mcL \otimes \mcN)
\end{align}
forms a dormant $\mr{GL}_2$-oper on $X$.
Then, we can define  the equivalence relation ``$\sim$" in the set of dormant $\mr{GL}_2$-opers on $X$ defined by the condition that  $\mcF^\heartsuit \sim \mcF'^\heartsuit$ if and only if $\mcF^\heartsuit_{\otimes \msN}  \cong \mcF'^\heartsuit$ for some   $\msN$ as above. 

 Let us write
   \begin{align}
   \mr{Op}_{\mfs \mfl_2, X}^{^\mr{Zzz...}} \  \left(\text{resp.,} \ \mr{O}\mr{p}_{\mr{GL}_2, X, /\sim}^{^\mr{Zzz...}} \right)
   \end{align}
  for the set  of isomorphism classes of dormant $\mfs \mfl_2$-opers (resp.,  the set of equivalence classes of  dormant $\mr{GL}_2$-opers) on $X$.


Each (dormant) $\mr{GL}_2$-oper $\mcF^\heartsuit$ induces, via projectivization, a(n) (dormant) $\mfs \mfl_2$-oper; we shall denote the resulting $\mfs \mfl_2$-oper by  
 \begin{align}
 \mcF^{\heartsuit \Rightarrow \spadesuit}.
 \end{align}
According to ~\cite[Theorem D]{Wak5}, 
the assignment $\mcF^\heartsuit \mapsto \mcF^{\heartsuit \Rightarrow \spadesuit}$ defines a well-defined bijection  
\begin{align} \label{wE78}
\mr{O}\mr{p}_{\mr{GL}_2, X, /\sim}^{^\mr{Zzz...}} \isom \mr{Op}_{\mfs \mfl_2, X}^{^\mr{Zzz...}};
\end{align}
 this  means that, for a dormant 
 $\mfs \mfl_2$-oper $\mcE^\spadesuit$, there exists  a dormant $\mr{GL}_2$-oper $\mcF^\heartsuit$ with $\mcF^{\heartsuit \Rightarrow \spadesuit} \cong \mcE^\spadesuit$ and the equivalence class of  $\mcF^\heartsuit$
  is uniquely determined.

 \SSP
  \bde \label{wD3}
  \begin{itemize}
  \item[(i)]
A {\bf dormant  curve} over $k$  is a pair
\begin{align}
\mbX :=  (X, \mcE^\spadesuit)
\end{align}
consisting of a connected smooth proper  curve $X$ over $k$  and a dormant $\mfs \mfl_2$-oper $\mcE^\spadesuit$ on $X$.
\item[(ii)]
Let $\mbX := (X, \mcE_X^\spadesuit)$ and $\mbY := (Y, \mcE_Y^\spadesuit)$ be dormant curves over $k$.
An {\bf \'{e}tale covering} from $\mbY$ to $\mbX$ is an \'{e}tale covering $w : Y \migi X$ with  $w^* (\mcE_X^\spadesuit) \cong \mcE_Y^\spadesuit$ (cf. ~\cite[Chap.\,2, \S\,2.1.5]{Wak5} for the definition of the pull-back of an  $\mfs \mfl_2$-oper); for convenience, we denote such a morphism by  $w : \mbY \migi \mbX$.
Also, we shall say that an \'{e}tale covering $w : \mbY \migi \mbX$ is {\bf Galois} (resp., {\bf cyclic}) if $w : Y \migi X$ is Galois (resp., cyclic). 
If $w : \mbY \migi \mbX$ is a Galois \'{e}tale covering, then we shall write $\mr{Gal}(\mbY/\mbX) := \mr{Gal} (Y/X)$. 

\end{itemize}
\ede
 
 \SSP
  \bde \label{wD2}
\begin{itemize}
\item[(i)]
Let $\mbX := (X, \mcE^\spadesuit)$ (where $\mcE^\spadesuit := (\mcE, \nabla, \mcE_B)$) be a dormant curve over $k$.
We shall say that $\mbX$ (or $\mcE^\spadesuit$) is {\bf ordinary} if the morphism 
\begin{align} \label{wE102}
\Xi_{\mcE^\spadesuit} : H^1 (X, \mr{Ker}(\nabla^\mr{ad})) \migi H^1 (X,  (\mfs \mfl_2/\mfb_2)_{\mcE_B})
\end{align}
 induced by the composite of natural morphisms $\mr{Ker}(\nabla^\mr{ad}) \migiincl (\mfs \mfl_2)_\mcE \left(=(\mfs \mfl_2)_{\mcE_B} \right) \migisurj (\mfs \mfl_2/\mfb_2)_{\mcE_B}$ is an isomorphism.
 \item[(ii)]
 Let $X$ be a  connected smooth proper curve over $k$.
 We shall say that $X$ is {\bf dormant-ordinary}
 if any dormant $\mfs \mfl_2$-oper on $X$ is ordinary.
\end{itemize}
\ede

\SSP

\begin{rema}\label{wR100}
Let 
$\mcE^\spadesuit := (\mcE, \nabla, \mcE_B)$ be a dormant  $\mfs \mfl_2$-oper on $X$,  and choose a dormant 
 $\mr{GL}_2$-oper $\mcF^\heartsuit := (\mcF, \nabla_\mcF, \mcL)$  with $\mcF^{\heartsuit \Rightarrow \spadesuit} \cong \mcE^\spadesuit$.
The sheaf of horizontal sections $\mcF^\nabla$ in $\mcF$ with respect to $\nabla_\mcF$ (i.e., $\mcF^\nabla := \mr{Ker}(\nabla_\mcF)$) may be regarded as an $\mcO_{X^{(1)}}$-module via the underlying homeomorphism of $F_{X/k}$.
 Denote by $\mcE nd^0 (\mcF)$ (resp., $\mcE nd^0 (\mcF^\nabla)$) the sheaf of $\mcO_{X}$-linear  (resp., $\mcO_{X^{(1)}}$-linear) endomorphisms of $\mcF$ (resp., $\mcF^\nabla$)  with vanishing trace.
 One may verify that  $\mr{Ker}(\nabla^\mr{ad})$ and $(\mfs \mfl_2/\mfb_2)_{\mcE_B}$ are  naturally identified with  $\mcE nd^0 (\mcF^\nabla)$ and $\mcH om_{\mcO_X}(\mcL, \mcF/\mcL)$ respectively. 
 Under these identifications,
   the morphism $\Xi_{\mcE^\spadesuit}$ coincides with  the morphism
 \begin{align}\label{wE101}
 H^1 (X^{(1)}, \mcE nd^0 (\mcF^\nabla)) \migi H^1 (X, \mcH om_{\mcO_X}(\mcL, \mcF/\mcL))
 \end{align}
 induced by the composite
 \begin{align} \label{wE100}
 F_{X/k}^{-1}(\mcE nd^0 (\mcF^\nabla)) \migiincl  \mcE nd^0 (\mcF) \migisurj  \mcH om_{\mcO_X}(\mcL, \mcF/\mcL),
 \end{align}
 where
 the first arrow arises from pull-back by $F_{X/k}$ (under the natural identification $F^*_{X/k}(\mcF^\nabla) \isom \mcF$) and 
  the second arrow arises from  both $\mcL \migiincl \mcF$ and  $\mcF \migisurj \mcF/\mcL$.

\end{rema}

\SSP

\begin{rema} \label{wR28}
Let  us fix an integer $g>1$, and denote by
$\mfM_g$ the moduli stack classifying connected smooth proper curves over $k$ of genus $g$.
Also, denote by $\mfM_g^{^\mr{Zzz...}}\!$ the moduli stack classifying dormant curves over $k$ of genus $g$.
Note   that a dormant curve $\mbX := (X, \mcE^\spadesuit)$  is ordinary if and only if  
$\mfM_g^{^\mr{Zzz...}}\!$ is \'{e}tale over $\mfM_g$ at the  point classifying $\mbX$  (cf.  ~\cite[\S\,6, Proposition  6.7.2]{Wak11}).
That is to say,  if $\mr{Def}_{\mbX}$ (resp., $\mr{Def}_X$) denotes the space of first order deformations of 
 $\mbX$ (resp., $X$), then the ordinariness of $\mbX$ is equivalent to the bijectivity of   the map 
 \begin{align}\label{wE57}
 \mr{Def}_{\mbX} \migi \mr{Def}_X
 \end{align}
 given by forgetting the data of the deformations of $\mcE^\spadesuit$.
 (This equivalence remains true even when $g=0$, $1$.)

In the context of $p$-adic Teichm\"{u}ller theory, dormant $\mfs \mfl_2$-opers were substantially investigated by S. Mochizuki (cf. ~\cite{Mzk2})  under the name of dormant indigenous bundles.
In that work,  he proved that $\mfM_g^{^\mr{Zzz...}}\!$ is a connected smooth Deligne-Mumford stack  over $k$  and the natural projection $\mfM_g^{^\mr{Zzz...}} \!\migi \mfM_g$ is finite, faithfully flat, and generically \'{e}tale (cf. ~\cite[Chap.\,II, \S\,2.3, Theorem 2.8]{Mzk2}).
In particular, a general dormant curve of genus $>1$ is ordinary, and moreover,  a general curve of genus $>1$ is dormant-ordinary.

One may also verify that the complement of the ordinary locus  forms a divisor on $\mfM_g^{^\mr{Zzz...}}\!$.
Indeed,  if  $(f: \widetilde{X}\migi \mfM_g^{^\mr{Zzz...}}\!, \widetilde{\mcE}^\spadesuit)$ (where $\widetilde{\mcE}^\spadesuit := (\widetilde{\mcE}, \widetilde{\nabla}, \widetilde{\mcE}_B)$)  denotes the universal object over $\mfM^{^\mr{Zzz...}}_g\!$,
then we obtain the universal  family
\begin{align}
\Xi_{\widetilde{\mcE}^\spadesuit} : \mbR^1 f_* (\mr{Ker}(\widetilde{\nabla}^\mr{ad})) \migi \mbR^1 f_* ((\mfs \mfl_2/\mfb_2)_{\widetilde{\mcE}_B})
\end{align}
 of ``$\Xi_{\mcE^\spadesuit}$"; this morphism  produces 
 a  morphism 
  \begin{align}
  \bigwedge^{3g-3}\Xi_{\widetilde{\mcE}^\spadesuit} : \bigwedge^{3g-3} \mbR^1 f_* (\mr{Ker}(\widetilde{\nabla}^\mr{ad}))  \migi \bigwedge^{3g-3} \mbR^1 f_* ((\mfs \mfl_2/\mfb_2)_{\widetilde{\mcE}_B})
  \end{align}
between line bundles on $\mfM_g^{^\mr{Zzz...}}\!$, which corresponds to 
 a global section  of the line bundle 
 \begin{align}
 \left(\bigwedge^{3g-3} \mbR^1 f_* ((\mfs \mfl_2/\mfb_2)_{\widetilde{\mcE}_B}) \right) \otimes \left(\bigwedge^{3g-3} \mbR^1 f_* (\mr{Ker}(\widetilde{\nabla}^\mr{ad})) \right)^\vee.
 \end{align}
 The morphism  $\Xi_{\widetilde{\mcE}^\spadesuit}$ is an isomorphism exactly on the complement of the divisor $\mfD \left(\subseteq \mfM_g^{^\mr{Zzz...}}\! \right)$ associated to this section.
 Thus,  the nonordinary locus in $\mfM_g^{^\mr{Zzz...}}\!$ coincides with the support  of $\mfD$. 
\end{rema}

\SSP

\begin{rema} \label{wR29}
In the study of $p$-adic Teichm\"{u}ller theory,
there is   an important class of indigenous bundles, i.e.,
{\it (hyperbolically)  ordinary} nilpotent indigenous bundles (cf. ~\cite[Chap.\,II, \S\,3, Definition 3.1]{Mzk1}).
Despite the similarity in concept,
it has nothing to do with the ordinariness discussed in the present paper.
Indeed,  for an integer $g>1$, let  $\mfN_{g}$ (resp., ${^\circledcirc}\mfN_g$) denote the moduli stack classifying
connected smooth proper  genus-$g$ curves over $k$  equipped with a nilpotent indigenous bundle (resp., an ordinary nilpotent indigenous bundle).
The projection $\mfN_g \migi \mfM_g$ is finite and faithfully flat of degree $p^{3g-3}$ (cf.  ~\cite[Chap.\,II, \S\,1.3, Proposition 1.7]{Mzk2}) and 
the substack ${^\circledcirc}\mfN_{g}$ of $\mfN_g$ coincides with the \'{e}tale locus relative to this projection.
Also, the stack $\mfM_g^{^\mr{Zzz...}}\!$
(cf. Remark \ref{wR28} above)
 may be considered as a closed substack of $\mfN_g$.
Since the $p$-curvature of any ordinary nilpotent indigenous  bundle does not vanish,
we have the equality $\mfM_g^{^\mr{Zzz...}}\! \cap {^\circledcirc}\mfN_g = \emptyset$ in $\mfN_g$.
\end{rema}


\SSP 
 \begin{exa}[Dormant $\mfs \mfl_2$-opers on the projective line]
 Recall  from ~\cite[\S\,3.5, Proposition 3.7]{Wak10} that there exists a unique dormant $\mfs \mfl_2$-oper on the projective line $\mbP := \mr{Proj}(k[x_1, x_2])$ over $k$;
 this is obtained as the projectivization $\mcF_0^{\heartsuit \Rightarrow \spadesuit}$ of the dormant $\mr{GL}_2$-oper 
 \begin{align}
 \mcF_0^\heartsuit := (\mcO_{\mbP}^{\oplus 2}, \nabla_0^{\oplus 2}, \iota_0),
 \end{align}
 where 
   $\iota_0$ denotes the $\mcO_{\mbP}$-linear injection $\mcO_{\mbP}(-1) \migiincl \mcO_{\mbP}^{\oplus 2}$ given by $v \mapsto (v \cdot  x_1, v \cdot  x_2)$ for any local section $v \in \mcO_{\mbP}(-1)$.
If we write $(\mcE, \nabla, \mcE_B) := \mcF_0^{\heartsuit \Rightarrow \spadesuit}$, then 
 it is verified  that $\mr{Ker}(\nabla^\mr{ad}) \cong \mcO_{\mbP^{(1)}}^{\oplus 3}$  and $(\mfs \mfl_2/\mfb_2)_{\mcE_B} \cong \mcO_{\mbP}(2)$, which implies  $H^1 (\mbP, \mr{Ker}(\nabla^\mr{ad})) = H^1 (\mbP, (\mfs \mfl_2/\mfb_2)_{\mcE_B}) =0$.
 In particular, the dormant  $\mfs \mfl_2$-oper $\mcF_0^{\heartsuit \Rightarrow \spadesuit}$ is ordinary.
 \end{exa}

\LSP
\subsection{Dormant $\mfs \mfl_2$-opers on an ordinary elliptic curve} \label{wS1}

In the rest of this section, we describe all possible dormant $\mfs \mfl_2$-opers of elliptic (i.e., genus-$1$) curves and classify ordinary  dormant elliptic curves.

First, let us consider the case where the underlying curve is ordinary.
Let $X$ be a connected smooth proper  genus-$1$ curve over $k$, and suppose that $X$ is ordinary, i.e., the morphism $F_{X/k}^\natural : H^1 (X^{(1)}, \mcO_{X^{(1)}}) \migi H^1 (X, \mcO_X)$ induced by $F_{X/k}$ is bijective.
Denote by
$\mr{Ver}_X^\times$ the set of isomorphism classes of line bundles $\mcL$ on $X^{(1)}$ with $F^*_{X/k}(\mcL) \cong \mcO_X$ and $\mcL \ncong \mcO_{X^{(1)}}$.
For each  line bundle  $\mcL$ classified by this set,
 $\nabla_\mcL^\mr{can}$ (cf. (\ref{wE44})) may be regarded as a $k$-connection on 
 $\mcO_X$  via an isomorphism  $F^*_{X/k}(\mcL) \isom \mcO_X$.
 Note that the resulting  $k$-connection on $\mcO_X$ does not depend on the choice  of  this isomorphism 
because of the equality $\Gamma (X, \mcO_X^\times) = k^\times \left(\subseteq \Gamma (X, F_{X/k}^{-1}(\mcO_{X^{(1)}}^\times)) \right)$. 
Also, denote by $\overline{\mr{Ver}}_X^\times$ the quotient set of $\mr{Ver}_X^\times$ by the equivalence class ``$\sim$" determined   by
$\mcL \sim \mcL^\vee$; given  each line bundle $\mcL$ classified by $\mr{Ver}_X^\times$,
we shall write $[\mcL]$ for the element of $\overline{\mr{Ver}}_X^\times$ represented by $\mcL$.
In particular,  the ordinariness of $X$ implies  the equality   $\sharp (\overline{\mr{Ver}}_X^\times) = \frac{p-1}{2}$.

For  each  line bundle  $\mcL$ classified by   
$\mr{Ver}_X^\times$,
we obtain the collection
\begin{align} \label{wE109}
\mcF_\mcL^\heartsuit := (\mcO_{X}^{\oplus 2}, \nabla_0 \oplus \nabla_\mcL^\mr{can}, \mr{Im}(\Delta)),
\end{align}
where  $\Delta$  denotes the diagonal embedding $\Delta : \mcO_X \migiincl \mcO_X^{\oplus 2}$.

\SSP
\ble \label{wL1}
\begin{itemize}
\item[(i)]
The collection  $\mcF_\mcL^\heartsuit$
 forms a dormant $\mr{GL}_2$-oper 
  on $X$.
Also,  if $\mcL$ and $\mcL'$ are line bundles classified by  $\mr{Ver}_X^\times$,
then $\mcF_\mcL^{\heartsuit \Rightarrow \spadesuit} \cong \mcF_{\mcL'}^{\heartsuit \Rightarrow \spadesuit}$ if and only if  $[\mcL] = [\mcL']$ in $\overline{\mr{Ver}}_X^\times$.

\item[(ii)]
For each $\mcL \in \mr{Ver}_X^\times$, the dormant $\mfs \mfl_2$-oper 
$\mcF_\mcL^{\heartsuit \Rightarrow \spadesuit}$ is ordinary.
\end{itemize}
\ele
\begin{proof}
First, let us consider  the former assertion of (i).
Since $\nabla_\mcL^\mr{can}$ has vanishing $p$-curvature,  the remaining portion is to 
 prove that $\mr{KS}_{\mcF_\mcL^\heartsuit} :  \mr{Im}(\Delta) \migi  \Omega_{X/k} \otimes \mr{Coker}(\Delta)$ is an isomorphism.
Suppose, on the contrary,  that $\mr{KS}_{\mcF_\mcL^\heartsuit}$ is not an isomorphism.
By the equality  $\mr{deg} (\mr{Im}(\Delta)) = \mr{deg}(\Omega_{X/k} \otimes \mr{Coker}(\Delta)) \left( =0\right)$,
this morphism must be the zero map.
This implies that the line subbundle $\mr{Im}(\Delta)$ of $\mcO_X^{\oplus 2}$ is closed under $\nabla_0 \oplus \nabla^\mr{can}_\mcL$.
But,  this is a contradiction because  the assumption 
$\mcL \ncong \mcO_{X^{(1)}}$ implies that
   $\nabla_0 \oplus \nabla^\mr{can}_\mcL$ sends $\begin{pmatrix} 1 \\ 1 \end{pmatrix}$ to $\partial \otimes \begin{pmatrix} 0 \\ a \end{pmatrix}$ (where  $\partial$ denotes  a generator of $\Omega_{X/k}$) for some $a \in k^\times$.
Hence,   $\mr{KS}_{\mcF_\mcL^\heartsuit}$ turns out to be an isomorphism, as desired.

Next, we shall consider  the latter assertion of (i).
The ``if'' part  follows immediately from (\ref{wE78}) and  the observation that 
the automorphism of $\mcO_X^{\oplus 2}$ given by switching the factors defines 
 an isomorphism  of $\mr{GL}_2$-opers  $\mcF^\heartsuit_\mcL \isom \mcF^\heartsuit_{\mcL^\vee, \otimes (\mcO_X, \nabla_\mcL^\mr{can})}$.
 
To prove the ``only if" part, we suppose that
$\mcF_{\mcL}^{\heartsuit \Rightarrow \spadesuit} \cong \mcF_{\mcL'}^{\heartsuit \Rightarrow \spadesuit}$.
By the bijectivity of (\ref{wE78}),
there exists 
an isomorphism  of $\mr{GL}_2$-opers  $\eta : \mcF^\heartsuit_\mcL \isom \mcF^\heartsuit_{\mcL', \otimes (\mcN, \nabla_\mcN)}$ for some line bundle  $\mcN$    on $X$  equipped with a $k$-connection $\nabla_\mcN$ with vanishing $p$-curvature.
This isomorphism is restricted to an isomorphism $\mr{Im}(\Delta) \isom \mr{Im}(\Delta) \otimes \mcN$  between the respective line subbundles.
It follows that $\mcN$ may be identified with $\mcO_X$, and that $\nabla_\mcN = \nabla_{\mcL''}^\mr{can}$ for some line bundle $\mcL''$ on $X^{(1)}$ with $F_{X/k}^*(\mcL'') \cong \mcO_X$.
The isomorphism $\eta$ is, moreover, restricted to an isomorphism $\mcO_{X^{(1)}}\oplus \mcL \isom \left((\mcO_{X^{(1)}} \oplus \mcL')\otimes \mcL'' = \right) \mcL'' \oplus \mcL' \otimes \mcL''$ between the respective sheaves of horizontal sections.
Hence, $\mcO_{X^{(1)}}$ is isomorphic to either $\mcL''$ or $\mcL' \otimes \mcL''$; the  former case implies  $\mcL \cong  \mcL'$ and the latter case implies $\mcL^\vee \cong \mcL''$.
At any rate, we have $[\mcL] = [\mcL']$ in $\overline{\mr{Ver}}^\times_X$, as desired.

Next, we shall consider  assertion (ii).
For simplicity, we write 
$\mcF := \mcO_X^{\oplus 2}$,  $\nabla := \nabla_0 \oplus \nabla_\mcL^\mr{can}$.
Since $\mr{Im}(\Delta) \cong \mr{Coker}(\Delta) \cong  \mcO_X$,
we have 
\begin{align} \label{wE105}
H^1 (X, \mcH om_{\mcO_X}(\mr{Im}(\Delta), \mr{Coker}(\Delta))) \cong H^1 (X, \mcO_X).
\end{align}
On the other hand, the sheaf of horizontal sections $\mcF^\nabla \left(\subseteq \mcF \right)$ with respect to $\nabla$ coincides with
$\mcO_{X^{(1)}} \oplus \mcL$.
It follows that $\mcE nd^0 (\mcF^\nabla) \cong \mcO_{X^{(1)}} \oplus \mcL \oplus \mcL^\vee$.
The assumption $\mcL \ncong \mcO_{X^{(1)}}$ implies 
$H^1 (X, \mcL) = H^1 (X, \mcL^\vee) = 0$, so we have
\begin{align} \label{wE106}
H^1 (X^{(1)}, \mcE nd^0 (\mcF^\nabla)) \cong H^1 (X^{(1)}, \mcO_{X^{(1)}} \oplus \mcL \oplus \mcL^\vee) \cong H^1 (X^{(1)}, \mcO_{X^{(1)}}).
\end{align}
Under the identifications given by (\ref{wE105}) and (\ref{wE106}),
 $\Xi_{\mcF_\mcL^{\heartsuit \Rightarrow \spadesuit}}$ coincides with the morphism $F_{X/k}^\natural$ up to multiplication by an element of $k^\times$.
In particular, by the ordinariness of $X$, $\Xi_{\mcF_\mcL^{\heartsuit \Rightarrow \spadesuit}}$ turns out to be an isomorphism, i.e., 
 $\mcF^{\heartsuit \Rightarrow \spadesuit}_\mcL$
  is ordinary.
This completes the proof of this lemma.
\end{proof}

\SSP

\bpr \label{wP1}
Any dormant $\mfs \mfl_2$-oper is isomorphic to $\mcF_\mcL^{\heartsuit \Rightarrow \spadesuit}$ for some 
$\mcL \in \mr{Ver}_X^\times$.
In particular, 
the assignment $\mcL \mapsto \mcF_\mcL^{\heartsuit \Rightarrow \spadesuit}$ induces  a bijection of sets
\begin{align}
\overline{\mr{Ver}}_X^{\times} \isom \mr{Op}_{\mfs \mfl_2, X}^{^\mr{Zzz...}}.
\end{align}
(hence $\sharp (\mr{Op}_{\mfs \mfl_2, X}^{^\mr{Zzz...}}) = \frac{p-1}{2}$).
\epr
\begin{proof}
Let $\mcE^\spadesuit$ be a dormant $\mfs \mfl_2$-oper on $X$.
We shall  apply  ~\cite[Chap.\,4, \S\,4.11, Corollary 4.70]{Wak5} in  the case where $n=2$ and  the $2$-theta characteristic ``$(\varTheta, \nabla_\vartheta)$" is taken as $(\mcO_X, \nabla_0)$.
(Note that the statement of {\it loc.\,cit}. remains true even when the underlying curve is not hyperbolic.)
Then, there exists  a dormant $\mr{GL}_2$-oper  $\mcF^\heartsuit := (\mcF, \nabla, \mcL)$  on $X$ such that $\mcF^{\heartsuit \Rightarrow \spadesuit} \cong \mcE^\spadesuit$, 
$\mcL = \mcF/\mcL = \mcO_X$, 
and the determinant of $(\mcF, \nabla)$ is isomorphic to $(\mcO_X, \nabla_0)$.
In particular,  the sheaf of horizontal sections $\mcF^\nabla \left(\subseteq \mcF \right)$    with respect to $\nabla$ satisfies  $\mr{det}(\mcF^\nabla) \cong \mcO_{X^{(1)}}$.
Moreover, according to the discussion in ~\cite[Chap.\,4, \S\,4.4, Example 4.20]{Wak5},
the extension $0 \migi \mcL  \migi \mcF \migi  \mcF/\mcL  \migi 0$
 splits, i.e., $\mcF \cong \mcL \oplus (\mcF/\mcL) \cong \mcO_X^{\oplus 2}$.

Here, let us suppose that $\mcF^\nabla$ is indecomposable.
By 
~\cite[\S\,10.2, Proposition 10.48]{Muk}, there exists a line bundle $\mcN$ on $X^{(1)}$ of degree $0$ which fits into an extension
$e_\mcN := (0 \migi \mcN \migi \mcF^\nabla \migi \mcN \migi 0)$.
The pull-back of $e_\mcN$ by $F_{X/k}$ defines an extension  $e^F_\mcN := (0 \migi F_{X/k}^*(\mcN) \migi \mcF \migi F_{X/k}^*(\mcN) \migi 0)$.
Since $\mr{KS}_{\mcF^\heartsuit}$
is an isomorphism, 
the composite $\mcL \migiincl \mcF \migisurj F_{X/k}^*(\mcN)$ is not the zero map.
The equality $\mr{deg}(\mcL) = \mr{deg}(F^*_{X/k}(\mcN)) \left(=0\right)$ implies that
this composite is an isomorphism.
Hence, $e^F_\mcN$ splits, which means 
 that   the element   $[e_\mcN]$ of $H^1 (X^{(1)}, \mcO_{X^{(1)}}) \left(= \mr{Ext}^1(\mcN, \mcN) \right)$ represented by $e_\mcN$ is mapped to  $0\in H^1 (X, \mcO_X)  \left(= \mr{Ext}^1(F_{X/k}^*(\mcN), F_{X/k}^*(\mcN)) \right)$ via  $F_{X/k}^\natural$.
But, the ordinariness of $X$ implies that  $[e_\mcN]$ must be the zero element, or equivalently, $e_\mcN$ splits.
This contradicts the assumption that $\mcF^\nabla$ is indecomposable.

By this fact  and $\mr{det}(\mcF^\nabla) \cong \mcO_{X^{(1)}}$, 
there exist a line bundle $\mcN'$ on $X^{(1)}$ and an isomorphism  $\eta : \mcF^\nabla \isom  \mcN'^\vee \oplus   \mcN'$.
Since $\mcF \left(= F_{X/k}^*(\mcF^\nabla) \right) \cong \mcO_X^{\oplus 2}$ as mentioned  above,
we have $F_{X/k}^*(\mcN') \cong \mcO_X$.
It follows that 
the pull-back $F_{X/k}^*(\eta)$ of $\eta$ defines an isomorphism
 $\mcF \isom \mcO_X^{\oplus 2}$ compatible with
the respective $k$-connections $\nabla$ and $\nabla_{\mcN'^\vee}^\mr{can} \oplus \nabla_{\mcN'}^\mr{can}$.
Since $\mr{KS}_{\mcF^\heartsuit}$ is an isomorphism,
the composite $\mcO_X \left(= \mcL \right) \migiincl \mcF \xrightarrow{F_{X/k}^*(\eta)} \mcO_X^{\oplus 2}$ does not  coincide with the inclusion into  the first nor second factor.
After possibly composing $F_{X/k}^*(\eta)$  with an automorphism of $\mcO_X^{\oplus 2}$,
we may assume that $F_{X/k}^*(\eta) (\mcL) = \mr{Im}(\Delta)$.
Hence, $F_{X/k}^*(\eta)$ defines an isomorphism $\mcF^\heartsuit_{\otimes (\mcO_X, \nabla_{\mcN'}^\mr{can})} \isom  \mcF_{\mcN'^{\otimes 2}}^\heartsuit$, which implies that $\mcF_{\mcN'^{\otimes 2}}^{\heartsuit \Rightarrow \spadesuit} \cong  \mcF^{\heartsuit \Rightarrow \spadesuit}_{\otimes (\mcO_X, \nabla_{\mcN'}^\mr{can})}  \cong \mcF^{\heartsuit \Rightarrow \spadesuit} \cong  \mcE^\spadesuit$.
This completes the proof of the former assertion.
The latter assertion follows immediately from the former assertion and Lemma \ref{wL1}, (i).
\end{proof}

\LSP
\subsection{Dormant $\mfs \mfl_2$-opers on a nonordinary elliptic curve} \label{wS3}

Next, let us consider the case where the underlying curve is nonordinary (= supersingular).
Let $X$ be a connected smooth proper genus-$1$ curve over $k$, and suppose that $X$ is nonordinary, i.e., $F_{X/k}^\natural =0$.

 Let $e := (0 \migi \mcO_{X^{(1)}}  \xrightarrow{\nu^\flat} \mcV \xrightarrow{\nu^\sharp} \mcO_{X^{(1)}} \migi 0)$ be a nontrivial extension of $\mcO_{X^{(1)}}$ by $\mcO_{X^{(1)}}$ itself; this is uniquely determined up to isomorphism since $\mr{Ext}^1 (\mcO_{X^{(1)}}, \mcO_{X^{(1)}}) \left(\cong H^1 (X^{(1)}, \mcO_{X^{(1)}})\right) = k$.
 This extension yields a $2$-step decreasing filtration $\{ \mcV^j \}_{j=0}^2$  of $\mcV$ given by $\mcV^0 = \mcV$, $\mcV^1 = \mcO_{X^{(1)}}$, and $\mcV^2 =0$.
 If $\mcE nd^0 (\mcV)$ denotes  the sheaf of $\mcO_{X^{(1)}}$-linear endomorphisms of $\mcV$ with vanishing trace,
 then the filtration $\{ \mcV^j\}_j$ induces a filtration $\{ \mcE nd^0 (\mcV)^j \}_j$ on this sheaf given by
 \begin{align}
\mcE nd^0(\mcV)^j := \{ f \in \mcE nd^0 (\mcF) \, | \, f (\mcV^l) \subseteq \mcV^{l+j} \ \text{for any $l$} \}.
\end{align}
For each $j = -1, 0, 1$,
the subquotient $\mcE nd^0(\mcV)^j/\mcE nd^0(\mcV)^{j+1}$ is naturally isomorphic to $\mcO_{X^{(1)}}$.
In particular, we obtain the surjection
\begin{align} \label{wE199}
\delta : \mcE nd^0(\mcV) \migisurj \left(\mcE nd^0(\mcV)^{-1}/\mcE nd^0(\mcV)^{0} = \right) \mcO_{X^{(1)}}.
\end{align}
Also,  the assignment from  each local section $f \in \mcE nd^0(\mcV)$ to  the  image  $f(1)$ of $1 \in \mcO_{X^{(1)}} \subseteq \mcV$ via this morphism 
 defines 
\begin{align} \label{wE107}
\mcE nd^0(\mcV) \migi \mcV.
\end{align}
Here, let us observe the following lemma, which will be used in the proof of  Lemma \ref{wL2}.

 \SSP
\ble \label{wL4}
The morphism $H^1 (\delta) : H^1 (X^{(1)}, \mcE nd^0(\mcV)) \migi H^1 (X^{(1)}, \mcO_{X^{(1)}})$ induced by $\delta$ is an isomorphism.
\ele
\begin{proof}
The extension $e$ induces the long exact sequence 
\begin{align}
\Gamma (X^{(1)}, \mcO_{X^{(1)}}) \migi H^1 (X, \mcO_{X^{(1)}}) \migi H^1 (X^{(1)}, \mcV) \migi H^1 (X^{(1)}, \mcO_{X^{(1)}}) \migi 0.
\end{align}
The first arrow sends $1$ to the class $[e] \left(\neq 0 \right)$ represented by $e$.
It follows that the third arrow is an isomorphism.
Next, 
for each local section $v \in \mcV$, we shall write $f_v$ for the locally defined endomorphism of $\mcV$ given by $a \mapsto a \cdot \overline{v} -  \overline{a} \cdot v$, where $\overline{(-)}$ denotes the natural projection $\mcV \migisurj \mcO_{X^{(1)}}$.
Then, the assignment $v \mapsto f_v$ defines an $\mcO_{X^{(1)}}$-linear injection $\mcV \migiincl \mcE nd^0 (\mcV)$, which  induces  an isomorphism $\mcV \isom \mcE nd^0 (\mcV)^0$.
This injection
fits into the following morphism of short exact sequences:
\begin{align}
\vcenter{\xymatrix@C=46pt@R=36pt{
0 \ar[r] & \mcV \ar[r]^-{v \mapsto f_v} \ar[d]^-{\mr{projection}} & \mcE nd^0 (\mcV) \ar[r]^-{\delta} \ar[d]^-{(\ref{wE107})} & \mcO_{X^{(1)}}\ar[r] \ar[d]_-{\wr}^-{\mr{id}} & 0
\\
0 \ar[r] & \mcO_{X^{(1)}}\ar[r]_-{\nu^\flat} & \mcV \ar[r]_-{\nu^\sharp} & \mcO_{X^{(1)}} \ar[r]& 0.
}}
\end{align}
This diagram induces a morphism between the associated long exact sequences:
\begin{align}
\vcenter{\xymatrix@C=16pt@R=36pt{
\Gamma (X^{(1)}, \mcO_{X^{(1)}})\ar[r]  \ar[d]^-{\mr{id}}_-{\wr} & H^1(X^{(1)}, \mcV) \ar[r] \ar[d] & H^1(X^{(1)}, \mcE nd^0 (\mcV)) \ar[r]^-{H^1 (\delta)} \ar[d] & H^1(X^{(1)}, \mcO_{X^{(1)}})\ar[r] \ar[d]_-{\wr}^-{\mr{id}} & 0
\\
\Gamma (X^{(1)}, \mcO_{X^{(1)}}) \ar[r] & H^1(X^{(1)}, \mcO_{X^{(1)}})\ar[r] & H^1(X^{(1)}, \mcV)  \ar[r]_-{\sim} & H^1(X^{(1)}, \mcO_{X^{(1)}}) \ar[r]& 0.
}}
\end{align}
The above observation shows that both the leftmost lower horizontal arrow and the second vertical arrow from the left are isomorphisms.
Hence,  the leftmost upper horizontal arrow becomes an isomorphism,  and $H^1 (\delta)$ turns out to be an isomorphism.
\end{proof}
 \SSP
 
  Note that the pull-back $e^F$ of $e$ by $F_{X/k}$ splits because of the assumption $F_{X/k}^\natural =0$.
  The set of split injections $\mcO_X \migiincl F_{X/k}^*(\mcV)$ of $e^F$  admits naturally a structure of  torsor modeled on 
 the $k$-vector space  $\mr{Hom}_{\mcO_X}(\mcO_X, \mcO_X) = k$.
By choosing one of such injections $\iota$, 
 we obtain the collection
   \begin{align}
 \mcF_{\otimes}^\heartsuit :=  (F^*(\mcV), \nabla_\mcV^\mr{can}, \mr{Im}(\iota)).
 \end{align}
 
 \SSP
\ble \label{wL2}
\begin{itemize}
\item[(i)]
The collection $\mcF_{\otimes}^\heartsuit$ forms a dormant $\mr{GL}_2$-oper on $X$.
Moreover, the isomorphism class of $\mcF_\otimes^\heartsuit$ does not depend on the choice of $\iota$.
\item[(ii)]
The dormant $\mfs \mfl_2$-oper $\mcF_\otimes^{\heartsuit \Rightarrow \spadesuit}$ is nonordinary.
\end{itemize}
\ele
\begin{proof}
First, let us consider the former assertion of (i).
Since $\nabla_\mcV^\mr{can}$ has vanishing $p$-curvature,
it suffices to prove that
$\mr{KS}_{\mcF_\otimes^\heartsuit} : \mr{Im}(\iota) \migi \Omega_{X/k} \otimes \mr{Coker}(\iota)$ is an isomorphism.
Suppose, on the contrary, that $\mr{KS}_{\mcF_\otimes^\heartsuit}$ is not an isomorphism.
By the equality $\mr{deg}(\mr{Im}(\iota)) = \mr{deg}(\Omega_{X/k} \otimes \mr{Coker}(\Delta)) \left(= 0\right)$,
this morphism must be the zero map.
This implies that the line subbundle $\mr{Im}(\Delta)$ of $F_{X/k}^*(\mcV)$
is closed under $\nabla_\mcV^\mr{can}$.
Since the $k$-connection $\nabla_\mcV^\mr{can}$ restricted to $\mr{Im}(\Delta)$ has vanishing $p$-curvature, 
the natural inclusion $\mr{Im}(\Delta) \cap \mcV \migi \mcV$ may be identified with $\iota$ after pull-back by  $F_{X/k}$ and  specifies a split injection  of $e$.
This contradicts the assumption that $e$ is a nontrivial extension. 
Thus,  $\mr{KS}_{\mcF_\otimes^\heartsuit}$ must be an isomorphism, as desired.

Next, we shall consider the latter assertion of (i).
Let us choose another split injection $\iota' : \mcO_X \migiincl  F_{X/k}^*(\mcV)$ of $e^F$.
In what follows, we prove that the two  $\mr{GL}_2$-opers $(F^*_{X/k}(\mcV), \nabla_\mcV^\mr{can}, \iota)$, $(F^*_{X/k}(\mcV), \nabla_\mcV^\mr{can}, \iota')$
are isomorphic.
Denote by $\mr{End}^\otimes_{\mcO_X} (F_{X/k}^*(\mcV))$ the subspace of $\mr{End}_{\mcO_X} (F_{X/k}^*(\mcV))$ consisting of endomorphisms preserving the extension structure $e^F$, i.e.,  inducing   the identity morphism of $\mcO_X$ via $F^*_{X/k}(\nu^\flat)$ and $F^*_{X/k}(\nu^\sharp)$. 
Then, there exists an element $h \in \mr{End}^\otimes_{\mcO_X} (F_{X/k}^*(\mcV))$ with $\iota' = h \circ \iota$.
One may verify that
the morphism of $k$-vector spaces
\begin{align}
\lambda_{X^{(1)}} : \Gamma (X^{(1)}, \mcO_{X^{(1)}}) \oplus \mr{End}_{\mcO_{X^{(1)}}}(\mcO_{X^{(1)}}) \migi  \mr{End}_{\mcO_{X^{(1)}}}(\mcV) \\
\left(\text{resp.,} \  \lambda_X :  \Gamma (X, \mcO_X) \oplus \mr{End}_{\mcO_X} (\mcO_X) \migi  \mr{End}^\otimes_{\mcO_X} (F_{X/k}^*(\mcV)) \right) \hspace{4mm} \notag
\end{align}
given by $(a, b) \mapsto a \cdot \mr{id}_{\mcV} + \nu^\flat \circ  b \circ \nu^\sharp$ (resp., $(a, b) \mapsto a \cdot \mr{id}_{\mcV} + F_{X/k}^*(\nu^\flat) \circ  b \circ F_{X/k}^*(\nu^\sharp)$) 
is an isomorphism.
(Note the non-resp'd portion follows from the assumption that $\mcV$ is indecomposable.)
The following square diagram is commutative:
\begin{align}
\vcenter{\xymatrix@C=96pt@R=36pt{
\Gamma (X^{(1)}, \mcO_{X^{(1)}}) \oplus \mr{End}_{\mcO_{X^{(1)}}}(\mcO_{X^{(1)}})\ar[r]^-{(a, b) \mapsto (F_{X/k}^*(a), F^*_{X/k}(b))} \ar[d]^-{\wr}_-{\lambda_{X^{(1)}}} & \Gamma (X, \mcO_X) \oplus \mr{End}_{\mcO_X} (\mcO_X) \ar[d]^-{\lambda_X}_-{\wr}
\\
\mr{End}_{\mcO_{X^{(1)}}}(\mcV)\ar[r]_-{f \mapsto F_{X/k}^*(f)} &  \mr{End}^\otimes_{\mcO_X} (F_{X/k}^*(\mcV)).
}}
\end{align}
Thus, the lower horizontal arrow turns out to be an isomorphism because of the bijectivity of the upper horizontal arrow.
In particular, there exists an $\mcO_{X^{(1)}}$-linear endomorphism $f$ of $\mcV$ with $F_{X/k}^*(f) = h$, which implies that
$h$ preserves $\nabla_\mcV^\mr{can}$.
Hence, $h$ defines an isomorphism  of $\mr{GL}_2$-opers $(F^*_{X/k}(\mcV), \nabla_\mcV^\mr{can}, \iota) \isom (F^*_{X/k}(\mcV), \nabla_\mcV^\mr{can}, \iota')$.
This completes the proof of the latter assertion of (i).

Finally, we shall consider assertion (ii).
It follows from the various definitions  involved 
that the composite 
\begin{align}
\Xi_{\mcF_\otimes^{\heartsuit \Rightarrow \spadesuit}} \circ H^1(\delta)^{-1} : H^1 (X^{(1)}, \mcO_{X^{1}})\migi H^1 (X, \mcO_X) \left(= H^1 (X, \mcH om_{\mcO_X} (\mr{Im}(\iota), \mr{Coker}(\iota))) \right)
\end{align}
 (cf. Remark \ref{wR100}, Lemma \ref{wL4}) coincides with  $F_{X/k}^\natural$ up to multiplication by an element of $k^\times$.
But, since  $X$ is assumed to be nonordinary, this map must be the zero map.
It follows that $\mcF_\otimes^{\heartsuit \Rightarrow \spadesuit}$ is nonordinary.
This completes the proof of this assertion.
\end{proof}
 
 \SSP
\bpr \label{wP2}
Any dormant $\mfs \mfl_2$-oper is isomorphic to 
$\mcF_\otimes^{\heartsuit \Rightarrow \spadesuit}$, that is to say,  the following equality of sets holds:
\begin{align}
\mr{Op}_{\mfs \mfl_2, X}^{^\mr{Zzz...}} =  \left\{  \mcF_\otimes^{\heartsuit \Rightarrow \spadesuit} \right\}.
\end{align}
\epr
\begin{proof}
Let us take  a dormant $\mfs \mfl_2$-oper  $\mcE^\spadesuit$ on $X$.
We shall apply ~\cite[Chap.\,4, \S\,4.11, Corollary 4.70]{Wak5} in the case where $n=2$ and the $2$-theta characteristic ``$(\varTheta, \nabla_\vartheta)$" is taken as $(\mcO_X, \nabla_0)$.
(Note that the statement of {\it loc.\,cit}. remains true even when the underlying curve is not hyperbolic.)
Then, there exists a dormant $\mr{GL}_2$-oper $\mcF^\heartsuit := (\mcF, \nabla, \mcL)$ on $X$ such that $\mcF^{\heartsuit \Rightarrow \spadesuit} \cong \mcE^\spadesuit$, $\mcL = \mcF/\mcL = \mcO_X$, and the determinant of $(\mcF, \nabla)$ is isomorphic to $(\mcO_X, \nabla_0)$.
In particular, the sheaf of horizontal sections $\mcF^\nabla \left(\subseteq \mcF \right)$ with respect to $\nabla$ satisfies $\mr{det}(\mcF^\nabla) \cong \mcO_{X^{(1)}}$.
According to the discussion in ~\cite[Chap.\,4, \S\,4.4, Example 4.20]{Wak5},
the extension $0 \migi \mcL \migi \mcF \migi \mcF/ \mcL \migi 0$ splits, i.e., 
$\mcF \cong \mcL \oplus (\mcF/ \mcL) \cong \mcO_X^{\oplus 2}$.
Here, let us suppose that $\mcF^\nabla$ is decomposable.
Then, there exists a line bundle $\mcN$ on $X^{(1)}$  with $\mcF^\nabla \cong \mcN \oplus \mcN^\vee$.
Since $\mcO_X^{\oplus 2} \cong \mcF \cong F_{X/k}^*(\mcF^\nabla) \cong F_{X/k}^{*}(\mcN) \oplus F_{X/k}^{*}(\mcN^\vee)$,  
the pull-back $F^*_{X/k}(\mcN)$ of $\mcN$ must be isomorphic to  $\mcO_X$.
Hence, by the assumption that $X$ is nonordinary, $\mcN$  is isomorphic to $\mcO_{X^{(1)}}$,   so $\mcF^\nabla \cong \mcO_{X^{(1)}}^{\oplus 2}$.
This implies that $\nabla$ coincides with $\nabla_0^{\oplus 2}$ under a suitable  identification
$\mcF \isom \mcO_X^{\oplus 2}$.
In particular, any line subbundle of $\mcF$ is closed under $\nabla$.
This contradicts the assumption that $\mr{KS}_{\mcF^\heartsuit}$ is an isomorphism.
Thus, $\mcF^\nabla$ turns out to be  indecomposable,  in particular,   isomorphic to $\mcV$.
It follows that $\mcF^\heartsuit \cong \mcF_\otimes^\heartsuit$, which implies $\mcE^\spadesuit \cong \mcF_\otimes^{\heartsuit \Rightarrow \spadesuit}$.
This completes the proof of this proposition.
 \end{proof}
\SSP

Let us describe  a consequence of the discussions in \S\S\,\ref{wS1}-\ref{wS3}.

\SSP
 \bco
 Let $\mbX := (X, \mcE^\spadesuit)$ be a dormant smooth curve of genus $1$ over $k$.
 Then, $\mbX$ is ordinary if and only if $X$ is ordinary.
 In particular, any ordinary connected  smooth proper curve of genus $1$ is dormant-ordinary. 
 \eco

\vspace{10mm}
\section{Dormant $\mfs \mfl_n$-opers and Ordinariness} \SSP

This section 
deals with dormant $\mfs \mfl_n$-opers (for various  $n$) in the context of logarithmic geometry to formulate a generalized statement of Theorem \ref{EW1}, i.e., Theorem \ref{y019}.
We introduce 
 the ordinariness  of dormant $\mfs \mfl_n$-opers and describe   several  criteria for  that property.
By using one of these criteria, 
 we prove Theorem \ref{y019}, (i) (which includes  Theorem \ref{EW1}, (i)).
Let us fix an integer $n$ satisfying either  ``$n=2$ and $2<p$" or ``$2< 2n <p$".

\LSP
\subsection{Log curves and pointed stable curves} \label{S12}

Recall from ~\cite[Definition 4.5]{ACGH} that
a {\bf log curve} is a log smooth integrable morphism $f^\mr{log} : X^\mr{log} \migi S^\mr{log}$ of fs log schemes
 such that the geometric fibers of the underlying morphism $f : X \migi S$  are reduced connected $1$-dimensional schemes.

Denote by $\overline{\mfM}_{g,r}$  (cf. ~\cite[Definition 1.1]{Kn2})  the moduli stack of $r$-pointed stable curves over $k$ of genus $g$. 
 Also, denote  by $f_\mr{univ} : \mfC_{g,r} \migi \overline{\mfM}_{g,r}$ the universal family of  curves over $\overline{\mfM}_{g,r}$, which is equipped with  $r$ marked points $\sigma_{\mr{univ}, 1} , \cdots, \sigma_{\mr{univ}, r} : \overline{\mfM}_{g,r} \migi \mfC_{g,r}$.
It is well-known   that $\overline{\mfM}_{g,r}$ may be represented by a connected, proper, and smooth Deligne-Mumford stack over $k$ of dimension $3g-3+r$ (cf. ~\cite[\S\,5]{DM}, ~\cite[Corollary 2.6 and Theorem 2.7]{Kn2}).
The stack 
$\overline{\mfM}_{g,r}$ has a natural log structure given by the divisor at infinity (cf. ~\cite[Theorem 4.5]{KaFu}),
where we shall denote the resulting log stack  by $\overline{\mfM}_{g,r}^{\mr{log}}$.
Moreover, we obtain a log structure on $\mfC_{g,r}$ by taking   the divisor which is the union of the $\sigma_{\mr{univ}, i}$'s and the pull-back of the divisor at infinity of $\overline{\mfM}_{g,r}$ (resp., the divisor defined as the pull-back of the divisor at infinity of $\overline{\mfM}_{g,r}$); let us  denote  the resulting log stack by   $\mfC^{\mr{log}}_{g,r}$ (resp., $\mfC^{\mr{log}'}_{g,r}$).
Let $\mfM_{g,r}$ denotes the substack of $\overline{\mfM}_{g,r}$ classifying  smooth curves; it is a dense open substack of $\overline{\mfM}_{g,r}$ and coincides  with  the locus  in which the log structure of $\overline{\mfM}_{g,r}^\mr{log}$  becomes trivial.

Let us  fix a scheme $S$ over $k$ and
 an $r$-pointed stable curve
\begin{equation}
\label{X}
 \msX : =(f :X \migi S, \{ \sigma_i : S \migi X\}_{i=1}^r)\end{equation}
 of genus $g$ over $S$, which consists of a prestable curve $f : X \migi S$ over $S$ of genus $g$ and $r$ marked points $\sigma_i : S \migi X$ ($i = 1, \cdots, r$).
By pulling-back the  log structures of $\overline{\mfM}^{\mr{log}}_{g,r}$ and $\mfC^{\mr{log}}_{g,r}$ (resp.,  $\mfC^{\mr{log}'}_{g,r}$) via its classifying morphism,
we obtain log structures on $S$ and $X$ respectively; we denote  the resulting log schemes by 
\begin{align} \label{wE177}
S^{\msX\text{-}\mr{log}} \ \ \text{and}  \ \ X^{\msX\text{-}\mr{log}} \ \  \left(\text{resp.},  \ X^{\msX\text{-}\mr{log}'}\right).
\end{align} 
If there is no fear of confusion, we shall abbreviate 
them
 to $S^\mr{log}$ and $X^\mr{log}$ (resp., $X^{\mr{log}'}$) respectively.
The structure morphism $f : X \migi S$ extends to a morphism  of log schemes  $f^\mr{log} : X^{\mr{log}} \migi S^{\mr{log}}$ (resp., $f^{\mr{log}'} : X^{\mr{log}'} \migi S^{\mr{log}}$), 
by  which  $X^\mr{log}$ (resp., $X^{\mr{log}'}$)  becomes   a log curve over $S^\mr{log}$.
If, moreover, the underlying scheme $X$ is smooth over $S$, then $S^{\msX \text{-} \mr{log}} = S$.

An {\bf \'{e}tale (connected) covering} of $\msX$ (over $S$) is a pair $(\msY, w)$ consisting of   a pointed stable curve $\msY := (Y/S, \{ \sigma'_{j} \}_{j})$ over $S$ and an \'{e}tale covering $w : Y \migi X$ over $S$ with $w^{-1}(\bigcup_{i} \mr{Im}(\sigma_i)) = \bigcup_j \mr{Im}(\sigma_j)$. 
We shall say that an \'{e}tale covering $(\msY, w)$ is {\bf Galois} (resp., {\bf cyclic}) if $w$ is a Galois (resp., a cyclic) covering.
If the pair  $(\msY, w)$ as above  is a Galois covering of $\msX$ over $S$, then we shall write $\mr{Gal}(\msY/\msX) := \mr{Gal}(Y/X)$.

\LSP
\subsection{Dormant $\mfs \mfl_n$-opers and  dormant $\mr{GL}_n$-opers} \label{S21}

  Consider the Cartan decomposition $\mfs \mfl_n = \bigoplus_{j \in \mbZ} \mfs \mfl_{n, j}$ of $\mfs \mfl_n$, where $\mfb_n = \bigoplus_{j \in \mbZ_{\geq 0}} \mfs \mfl_{n, j}$;
 this decomposition gives rise to the decreasing  filtration $\{ \mfs \mfl_n^j \}_{j \in \mbZ}$ defined as $\mfs \mfl_n^j := \bigoplus_{l \geq j} \mfs \mfl_{n, l}$.
 Denote by $\Gamma$ the set of simple roots in $B_n$ with respect to the maximal torus of $\mr{PGL}_n$ consisting of elements represented by invertible diagonal matrices.
 Then,  there exists a natural decomposition 
 \begin{align} \label{wE113}
\mfs \mfl_n^{-1}/\mfs \mfl_n^0 = \bigoplus_{\alpha \in \Gamma} \mfs \mfl_n^{-\alpha},
 \end{align}
where $\mfs \mfl_n^{-\alpha}$ (for each $\alpha \in \Gamma$) denotes the root space of $- \alpha$.
Now, let    $S^\mr{log}$ be an fs  log scheme over $k$ and $f^\mr{log} : X^\mr{log} \migi S^\mr{log}$  a log curve over $S^\mr{log}$.

\vspace{3mm}
\bde \label{y0108}

 \begin{itemize}
 \item[(i)]
 Let us consider  a triple $\mcE^\spadesuit := (\mcE, \nabla, \mcE_B)$ consisting of a $\mr{PGL}_n$-bundle $\mcE$ on $X$, an $S^\mr{log}$-connection $\nabla$ on $\mcE$, and a $B_n$-reduction $\mcE_B$ of $\mcE$.
 We shall say that $\mcE^\spadesuit$ is an {\bf $\mfs \mfl_n$-oper} on $X^\mr{log}/S^\mr{log}$ if it satisfies the following two conditions:
 \begin{itemize}
 \item
 $\nabla (\mcT_{X^\mr{log}/S^\mr{log}}) \subseteq \widetilde{\mcT}^{-1}_{\mcE^\mr{log}/S^\mr{log}}$, where $\widetilde{\mcT}^{-1}_{\mcE^\mr{log}/S^\mr{log}} := (\mfs \mfl_n^{-1})_{\mcE_B} + \widetilde{\mcT}_{\mcE_B^\mr{log}/S^\mr{log}} \left(\subseteq \widetilde{\mcT}_{\mcE^\mr{log}/S^\mr{log}} \right)$;
 \item
 For each $\alpha \in \Gamma$, the composite
 \begin{align}
 \mcT_{X^\mr{log}/S^\mr{log}} \xrightarrow{\nabla} \widetilde{\mcT}^{-1}_{\mcE^\mr{log}/S^\mr{log}}  \migisurj \widetilde{\mcT}^{-1}_{\mcE^\mr{log}/S^\mr{log}} /\widetilde{\mcT}^0_{\mcE^\mr{log}/S^\mr{log}} \left(= (\mfs \mfl_n^{-1}/\mfs \mfl_n^0)_{\mcE_B} \right) \migisurj \mfg_{\mcE_B}^{-\alpha}
 \end{align}
 is an isomorphism, where the third arrow denotes the natural projection with respect to the decomposition (\ref{wE113}).
 \end{itemize}
 \item[(ii)]
 Let us consider a collection of data
 $\mcF^\heartsuit := (\mcF, \nabla, \{ \mcF^j \}_{j=0}^n)$, where
\begin{itemize}
\item
$\mcF$ is a rank $n$ vector bundle  on $X$;
\item
$\nabla$ is an $S^\mr{log}$-connection
 $\mcF \migi \Omega_{X^\mr{log}/S^\mr{log}} \otimes\mcF$ on $\mcF$;
\item
$\{ \mcF^j \}_{j=0}^n$ is   a decreasing filtration
 $0 = \mcF^n \subseteq \mcF^{n-1} \subseteq \dotsm \subseteq  \mcF^0= \mcF$
 on $\mcF$ by vector bundles.
\end{itemize}
We shall say that $\mcF^\heartsuit$ is  a {\bf $\mr{GL}_n$-oper} on  $X^\mr{log} /S^\mr{log}$ if it
 satisfies  the following three conditions:
\begin{itemize}
\item
The subquotients $\mcF^j / \mcF^{j+1}$ ($ 0\leq j\leq n-1$)  are line bundles;
\item
$\nabla (\mcF^j) \subseteq  \Omega_{X^\mr{log}/S^\mr{log}} \otimes \mcF^{j-1}$ ($1 \leq j \leq n-1$);
\item
The well-defined {\it $\mcO_X$-linear} morphisms
\begin{equation} \label{GL2}
\mr{KS}^j_{\mcF^\heartsuit} :  \mcF^j/\mcF^{j+1} \stackrel{}{\migi}   \Omega_{X^\mr{log}/S^\mr{log}} \otimes (\mcF^{j-1}/\mcF^j)\end{equation}
($1 \leq j \leq n-1$)
defined by assigning $\overline{a} \mapsto \overline{\nabla (a)}$ for any local section $a \in \mcF^j$ (where $\overline{(-)}$'s denote the images in the respective quotients) are isomorphisms.
\end{itemize}
\item[(iii)]
We shall say that an $\mfs \mfl_n$-oper $\mcE^\spadesuit := (\mcE, \nabla, \mcE_B)$ (resp., a $\mr{GL}_n$-oper $\mcF^\heartsuit := (\mcF, \nabla, \{ \mcF^j \}_{j=0}^n)$)
is {\bf dormant} if $\nabla$ has vanishing $p$-curvature.
\end{itemize}
  \ede
\SSP

If $X^\mr{log}/S^\mr{log}$ arises (in the sense of (\ref{wE177})) from a pointed stable curve $\msX$, then
we  refer to
any $\mfs \mfl_n$-oper (resp., $\mr{GL}_n$-oper) on $X^\mr{log}/S^\mr{log}$ as an {\it $\mfs \mfl_n$-oper on $\msX$} (resp., a {\it $\mr{GL}_n$-oper on $\msX$}).

Also. just as in the case of  rank $2$ opers on smooth curves,
we can define the notion of an isomorphism between two $\mfs \mfl_n$-opers (resp., $\mr{GL}_n$-opers).
Also, each (dormant) $\mr{GL}_n$-oper $\mcF^\heartsuit$ induces, via projectivization, a(n) (dormant) $\mfs \mfl_n$-oper; we shall denote the resulting $\mfs \mfl_n$-oper by $\mcF^{\heartsuit \Rightarrow\spadesuit}$.

\SSP

\begin{rema}[$n$-theta characteristics and $(\mr{GL}_n, \vartheta)$-opers] \label{wR110}
An {\bf $n$-theta characteristic} of $X^\mr{log}/S^\mr{log}$ (cf. ~\cite[Chap.\,4, \S\,4.6, Definition 4.31]{Wak5}) is a pair
$\vartheta := (\Theta, \nabla_\vartheta)$
consisting of a line bundle $\Theta$ on $X$ and an $S^\mr{log}$-connection $\nabla_\vartheta$ on the line bundle  $\mcT_{X^\mr{log}/S^\mr{log}}^{\otimes \frac{n (n-1)}{2}} \otimes \Theta^{\otimes n}$.
For example,  if we are given   a theta characteristic $\varTheta$ of  $X^\mr{log}/S^\mr{log}$ (i.e.,   a line bundle $\varTheta$  on $X$ equipped with an isomorphism $\iota_\varTheta : \mcO_X \isom \mcT_{X^\mr{log}/S^\mr{log}} \otimes \varTheta^{\otimes 2}$),
then the line bundle $\varTheta^{\otimes (n-1)}$ together with the $S^\mr{log}$-connection corresponding to $\nabla_0$ via $\iota_\varTheta^{\otimes \frac{n(n-1)}{2}}$ forms an $n$-theta characteristic.

Now, let us fix an $n$-theta characteristic $\vartheta := (\Theta, \nabla_\vartheta)$ of $X^\mr{log}/S^\mr{log}$.
For each integer $j$, we shall write $\mcD^{<j}$ for the sheaf of differential operators on $X^\mr{log}/S^\mr{log}$ of   order $< j$ (cf. ~\cite[Chap.\,4, \S\,4.2]{Wak5}).
Also,  write
\begin{align}
{^\dagger}\mcF_\Theta := \mcD^{<n} \otimes \Theta, \hspace{5mm}
{^\dagger}\mcF^j_{\Theta} := \mcD^{<n-j} \otimes \Theta \hspace{3mm} (j=0, 1, \cdots, n).
\end{align}
Then,  there exists  the natural  composite isomorphism
\begin{align} \label{wE123}
\mr{det}({^\dagger}\mcF_\Theta) \isom \bigotimes_{j=0}^{n-1} {^\dagger}\mcF_\Theta^j/{^\dagger}\mcF_\Theta^{j+1} \isom \bigotimes_{j=0}^{n-1} \mcT_{X^\mr{log}/S^\mr{log}}^{\otimes (n-1-j)} \otimes \Theta \isom \mcT_{X^\mr{log}/S^\mr{log}}^{\otimes \frac{n (n-1)}{2}} \otimes \Theta^{\otimes n}.
\end{align}

By a {\bf $(\mr{GL}_n, \vartheta)$-oper} on $X^\mr{log}/S^\mr{log}$ (cf. ~\cite[Chap.\,4, \S\,4.6, Definition 4.36]{Wak5}), we mean  a $\mr{GL}_n$-oper  of the form $({^\dagger}\mcF_\Theta, \nabla, \{ {^\dagger}\mcF_\Theta^j \}_j)$ satisfying  $\mr{det}(\nabla) = \nabla_\vartheta$ under the identification $\mr{det}({^\dagger}\mcF_\Theta) = \mcT_{X^\mr{log}/S^\mr{log}}^{\otimes \frac{n (n-1)}{2}} \otimes \Theta^{\otimes n}$ given by (\ref{wE123}).
If $\vartheta$ arises from a theta characteristic $\varTheta$ as above, then we shall refer to any  $(\mr{GL}_n, \vartheta)$-oper as a {\bf  $(\mr{GL}_n, \varTheta)$-oper}.

Denote by $\mr{Op}^{^\mr{Zzz...}}_{\mfs \mfl_n, X^\mr{log}/S^\mr{log}}$ (resp., $\mr{Op}^{^\mr{Zzz...}}_{\mr{GL}_n, \vartheta, X^\mr{log}/S^\mr{log}}$) the set of isomorphism classes of dormant $\mfs \mfl_n$-opers (resp., dormant $(\mr{GL}_n, \vartheta)$-opers)  on $X^\mr{log}/S^\mr{log}$.
It follows from ~\cite[Theorem D]{Wak5} that the assignment $\mcF^\heartsuit \mapsto \mcF^{\heartsuit \Rightarrow \spadesuit}$ defines a bijection of sets 
\begin{align} \label{wE134}
\mr{Op}^{^\mr{Zzz...}}_{\mr{GL}_n, \vartheta, X^\mr{log}/S^\mr{log}}  \isom \mr{Op}^{^\mr{Zzz...}}_{\mfs \mfl_n, X^\mr{log}/S^\mr{log}}.
\end{align}
\end{rema}

\LSP
\subsection{Pull-back of dormant $\mfs \mfl_n$-opers} \label{wSdS2}

Let $h^\mr{log} : Y^\mr{log} \migi T^\mr{log}$ be another log curve, and  suppose  that 
we are given a square diagram of log schemes
\begin{align} \label{ffg3Eg}
\vcenter{\xymatrix@C=46pt@R=36pt{
Y^\mr{log}\ar[r]^{w^\mr{log}} \ar[d]_-{h^\mr{log}} & X^\mr{log} \ar[d]^-{f^\mr{log}}
\\
T^\mr{log} \ar[r] & S^\mr{log}
}}
\end{align}
such that the induced morphism $Y^\mr{log} \migi X^\mr{log} \times_{S^\mr{log}}T^\mr{log}$ is log \'{e}tale.
(For example, if $(\msY, w)$ is an \'{e}tale covering of $\msX$ over $S$ and $Y$ denotes the underlying curve of $\msY$, then $w$ extends to  a log \'{e}tale morphism $w^\mr{log} : Y^\mr{log} \migi X^\mr{log}$.)
Also, let $\mcE^\spadesuit := (\mcE, \nabla, \mcE_B)$ be a dormant $\mfs \mfl_n$-oper on $X^\mr{log}/S^\mr{log}$.
The pull-back $w^*(\mcE)$  of $\mcE$ by $w$ defines a $\mr{PGL}_n$-bundle on $Y$, which admits the  $B_n$-reduction  $w^*(\mcE_B)$.
Moreover, since the morphism $w^*(\mcT_{X^\mr{log}/S^\mr{log}}) \migi \mcT_{Y^\mr{log}/T^\mr{log}}$ induced by $w^\mr{log}$ is an isomorphism,
$\nabla$ associates    
 a $T^\mr{log}$-connection  $w^*(\nabla)$  on $w^*(\mcE)$ via pull-back by $w$.
One may verify that the resulting collection
\begin{align}
w^*(\mcE^\spadesuit) := (w^*(\mcE), w^*(\nabla), w^*(\mcE_B))
\end{align}
forms  a dormant $\mfs \mfl_n$-oper on $Y^\mr{log}/T^\mr{log}$, which is called the {\bf pull-back} of $\mcE^\spadesuit$ by $w^\mr{log}$.
Similarly, we can construct the pull-back of a dormant $\mr{GL}_n$-oper by $w^\mr{log}$.

Next, suppose that $T = S$ and  (\ref{ffg3Eg}) is cartesian (hence $Y = X$).
Then, since $\mcT_{X^\mr{log}/S^\mr{log}} = \mcT_{Y^\mr{log}/T^\mr{log}}$,
 dormant $\mfs \mfl_n$-opers on $X^\mr{log}/S^\mr{log}$ may be regarded, via 
 pull-back  by $w^\mr{log}$, as  dormant $\mfs \mfl_n$-opers on $Y^\mr{log}/T^\mr{log}$.

\LSP
\subsection{Ordinariness of dormant $\mfs \mfl_n$-opers} \label{wSS2}

Let us fix a pair of nonnegative integers $(g,r)$ with $2g-2+r>0$.
Denote by $\overline{\mfO}  \mfp_{\mfs \mfl_n, g, r}^{^\mr{Zzz...}}$ (resp., $\mfO  \mfp_{\mfs \mfl_n, g, r}^{^\mr{Zzz...}}$) the moduli stack  classifying pairs $(\msX, \mcE^\spadesuit)$ consisting of an $r$-pointed stable (resp., smooth) curve $\msX$ and a dormant $\mfs \mfl_n$-oper $\mcE^\spadesuit$ on $\msX$.
According to ~\cite[Theorems C and  G]{Wak5}, $\overline{\mfO}  \mfp_{\mfs \mfl_n, g, r}^{^\mr{Zzz...}}$ may be represented by a proper Deligne-Mumford stack over $k$ and the projection $\overline{\mfO} \mfp_{\mfs \mfl_n, g,r}^{^\mr{Zzz...}} \migi \overline{\mfM}_{g,r}$ given by $(\msX, \mcE^\spadesuit) \mapsto \msX$ 
is finite and generically \'{e}tale (i.e., any irreducible component of $\overline{\mfO} \mfp_{\mfs \mfl_n, g,r}^{^\mr{Zzz...}}$ that dominates $\overline{\mfM}_{g,r}$ admits a dense open subscheme which is \'{e}tale over $\overline{\mfM}_{g,r}$).
If 
\begin{align}
{^\circledcirc}\overline{\mfO}  \mfp_{\mfs \mfl_n, g, r}^{^\mr{Zzz...}}
\end{align}
denotes the unramified locus in $\overline{\mfO}  \mfp_{\mfs \mfl_n, g, r}^{^\mr{Zzz...}}$ relative to $\overline{\mfM}_{g,r}$, then it coincides with the \'{e}tale locus  (cf. ~\cite[Chap.\,6, \S\,6.5,  Corollary 6.21]{Wak5}).
In particular, ${^\circledcirc}\overline{\mfO}  \mfp_{\mfs \mfl_n, g, r}^{^\mr{Zzz...}}$ forms an open substack of $\overline{\mfO} \mfp^{^\mr{Zzz...}}_{\mfs \mfl_n, g,r}$.

\SSP
\bde \label{wD043}
Let $\msX$ be a pointed stable curve over a $k$-scheme $S$ and $\mcE^\spadesuit$ a dormant $\mfs \mfl_n$-oper on $\msX$. 
We shall say that $\mcE^\spadesuit$
 is {\bf ordinary} if the image of the classifying morphism $S \migi \overline{\mfO} \mfp_{\mfs \mfl_n, g,r}^{^\mr{Zzz...}}$ lies in ${^\circledcirc}\overline{\mfO}  \mfp_{\mfs \mfl_n, g, r}^{^\mr{Zzz...}}$.
 (If $n=2$, then it follows from a comment in Remark \ref{wR28} that this notion is the same as the ordinariness defined in Definition \ref{wD2}, (i).)
 \ede
\SSP

\begin{rema} \label{wR134}
Let $\msX$ be a pointed stable curve over $k$ and $\mcE^\spadesuit$ a dormant $\mfs \mfl_n$-oper on $\msX$.
Then,  $\mcE^\spadesuit$ is ordinary if and only if for every first order deformation $\msX_\varepsilon$ of $\msX$, there exists a unique (up to isomorphism) deformation of $\mcE^\spadesuit$ over $\msX_\varepsilon$.  

\end{rema}
\SSP

\begin{rema} \label{wR135}
According to  ~\cite[Chap.\,8, \S\,8.5, Proposition  8.21]{Wak5},
any dormant $\mfs \mfl_n$-oper on  a $3$-pointed projective line (or more generally, a totally degenerate pointed curve) over $k$ is ordinary.
\end{rema}
\SSP

\SSP
\begin{rema}
\label{wR34}
Let us consider the relationship between the ordinariness of dormant $\mfs \mfl_2$-opers and that of dormant $\mfs \mfl_n$-opers.
After fixing an $\mfs \mfl_2$-triple in $\mfs \mfl_n$, we have a closed immersion $\overline{\mfO} \mfp^{^\mr{Zzz...}}_{\mfs \mfl_2, g,r} \migiincl \overline{\mfO} \mfp^{^\mr{Zzz...}}_{\mfs \mfl_n, g,r}$ (cf. ~\cite[Chap.\,2, \S\,2.6,  Theorem 2.24]{Wak5}).
In particular, this morphism  is restricted to an immersion ${^\circledcirc}\overline{\mfO} \mfp^{^\mr{Zzz...}}_{\mfs \mfl_2, g,r} \migiincl {^\circledcirc}\overline{\mfO} \mfp^{^\mr{Zzz...}}_{\mfs \mfl_n, g,r}$.
This means that if $\mcE^\spadesuit$ is a dormant $\mfs \mfl_2$-oper such that the associated $\mfs \mfl_n$-oper is ordinary, then $\mcE^\spadesuit$ itself is ordinary.
\end{rema}

\LSP
\subsection{Criteria for ordinariness } \label{wS9}

Let $\msX := (X/k, \{ \sigma_i \}_i)$  be a pointed stable curve over $k$ and $\mcE^\spadesuit := (\mcE, \nabla,  \mcE_B)$ a dormant $\mfs \mfl_n$-oper on $\msX$.
Denote by  
$\mr{Def}_{\mcE^\spadesuit}$ the space of first order deformations of $\mcE^\spadesuit$ (preserving the dormancy condition), i.e., the tangent space of 
$\overline{\mfO} \mfp^{^\mr{Zzz...}}_{\mfs \mfl_n, g,r}$ relative to $\overline{\mfM}_{g,r}$ at 
 the $k$-rational point classifying $(\msX, \mcE^\spadesuit)$. 

We shall denote by $\nabla^\mr{ad}$ the $k^\mr{log}$-connection on the adjoint bundle  $(\mfs \mfl_n)_{\mcE}$ induced by $\nabla$.
By the definition of an $\mfs \mfl_n$-oper, this connection  restricts  to a $k$-linear morphism 
\begin{align}
\nabla^{\mr{ad}(0)} : (\mfs \mfl_n)^0_{\mcE_B} \migi \Omega_{X^\mr{log}/k^\mr{log}} \otimes (\mfs \mfl_n)^{-1}_{\mcE_B}.
\end{align}
Write $\mcK^\bullet [\nabla^{\mr{ad}(0)}]$  for the complex of sheaves defined by $\nabla^{\mr{ad}(0)}$ concentrated at degrees $0$ and $1$.
Then, the composite of natural morphisms 
\begin{align}
\Omega_{X^\mr{log}/k^\mr{log}} \otimes (\mfs \mfl_n)^{-1}_{\mcE_B} \migiincl  \Omega_{X^\mr{log}/k^\mr{log}} \otimes (\mfs \mfl_n)_{\mcE_B} \left(= \Omega_{X^\mr{log}/k^\mr{log}} \otimes (\mfs \mfl_n)_{\mcE} \right) \migisurj \mr{Coker}(\nabla^\mr{ad})
\end{align}
 induces 
 a morphism of $k$-vector spaces
\begin{align} \label{wE33451}
\xi_{\mcE^\spadesuit} : \mbH^1 (X, \mcK^\bullet [\nabla^{\mr{ad}(0)}]) \migi \Gamma (X, \mr{Coker}(\nabla^\mr{ad})),
\end{align} 
where $\mbH^1 (X, \mcK^\bullet [\nabla^{\mr{ad}(0)}])$ 
denotes the first hypercohomology group of the complex $\mcK^\bullet [\nabla^{\mr{ad}(0)}]$.

\bpr \label{wP1231}
The $k$-vector space $\mr{Ker}(\xi_{\mcE^\spadesuit})$ is canonically isomorphic to
$\mr{Def}_{\mcE^\spadesuit}$.
In particular, $\mcE^\spadesuit$ is ordinary if and only if $\xi_{\mcE^\spadesuit}$ is injective.
\epr
\begin{proof}
The assertion follows from the discussion in ~\cite[Chap.\,6, \S\,6.5.2]{Wak5}.
\end{proof}
\SSP

In the following, we shall describe   several  criteria of ordinariness  under certain assumptions.

\SSP
\bpr \label{wP3}
Suppose that $\msX$ is unpointed (i.e., $\{ \sigma_i \}_i = \emptyset$) and $X$ is smooth over $k$.
Let us choose a theta characteristic $\varTheta$  on $X$.
Recall   that  there exists  a unique (up to isomorphism)
 dormant $(\mr{GL}_n, \varTheta)$-oper $\mcF^{\heartsuit} := (\mcF, \nabla_\mcF, \{ \mcF^j \}_{j})$  on $X$ with $\mcF^{\heartsuit \Rightarrow \spadesuit} \cong \mcE^\spadesuit$ (cf. Remark \ref{wR110}).
Also, denote by $\gamma_{\mcF^\heartsuit}$ the $\mcO_{X^{(1)}}$-linear composite
\begin{equation}
\gamma_{\mcF^\heartsuit} : \mr{Ker}(\nabla_\mcF) \xrightarrow{\mr{inclusion}}  F_{X/k*}(\mcF) \xrightarrow{\mr{quotient}} F_{X/k*} (\mcF/\mcF^1) \left(= F_{X/k*}(\varTheta^{\otimes (1-n)}) \right),
\end{equation}
which induces 
the following composite:
\begin{align} \label{wE23}
\delta_{\mcF^\heartsuit} : \mr{Hom}_{\mcO_{X^{(1)}}} (\mr{Ker}(\nabla_\mcF), \mr{Coker}(\gamma_{\mcF^\heartsuit}))\migi  H^1 (X^{(1)}, \mcE nd_{\mcO_{X^{(1)}}} (\mr{Ker}(\nabla_\mcF))) \migi H^1 (X^{(1)}, \mcO_{X^{(1)}}). \hspace{-8mm}
\end{align}
where
\begin{itemize}
\item
the first arrow denotes the connecting morphism in the long exact sequence associated to   the short exact sequence of $\mcO_{X^{(1)}}$-modules
\begin{align}
0 \migi \mcE nd_{\mcO_{X^{(1)}}} (\mr{Ker}(\nabla_\mcF)) \xrightarrow{f \mapsto \gamma_{\mcF^\heartsuit} \circ f} &\mcH om_{\mcO_{X^{(1)}}} (\mr{Ker}(\nabla_\mcF), F_{X/k*}(\mcF/\mcF^1)) \\
\migi & \mcH om_{\mcO_{X^{(1)}}} (\mr{Ker}(\nabla_\mcF), \mr{Coker}(\gamma_{\mcF^\heartsuit})) \migi 0;\notag
\end{align} 
\item
the  second arrow arises from   the  trace morphism  $\mcE nd_{\mcO_{X^{(1)}}} (\mr{Ker}(\nabla_\mcF)) \migi \mcO_{X^{(1)}}$.
\end{itemize}
Then, the $k$-vector space $\mr{Ker}(\gamma_{\mcF^\heartsuit})$ is canonically isomorphic to 
$\mr{Def}_{\mcE^\spadesuit}$.
In particular, $\mcE^\spadesuit$ is ordinary if and only if $\gamma_{\mcF^\heartsuit}$
is injective.
\epr
\begin{proof}
Denote by $\mr{Quot}^0$ (resp., $\mr{Quot}^\mr{triv}$) the Quot-scheme over $k$ classifying $\mcO_{X^{(1)}}$-submodules $\mcG$ of $F_{X/k*}(\varTheta^{\otimes (1-n)})$ with $\mr{rank} (\mcG) = 2$ and $\mr{deg}(\mcG)= 0$ (resp., $\mr{rank} (\mcG) = 2$ and  $\mr{det}(\mcG) \cong \mcO_{X^{(1)}}$). 
In particular,  $\mr{Quot}^\mr{triv}$ defines  a closed subscheme of $\mr{Quot}^0$.
According to  ~\cite[Chap.\,9, \S\,9.2, Proposition 9.4]{Wak5},  the assignment $\mcE^\spadesuit \mapsto \gamma_{\mcF^\heartsuit}$ defines an isomorphism of $k$-schemes
$\mfO \mfp_{\mfs \mfl_n, g, 0}^{^\mr{Zzz...}} \times_{\mfM_{g,0}, x} k \isom \mr{Quot}^{\mr{triv}}$, where $x$ denotes the $k$-rational point of $\mfM_{g, 0}$ classifying $X$.
In particular, $\gamma_{\mcF^\heartsuit}$ specifies a point of $\mr{Quot}^\mr{triv}$, as well as of $\mr{Quot}^{0}$.
It follows from  a well-known generalities on the deformation theory of Quot-schemes that the tangent space of $\mr{Quot}^0$ at the point classifying $\gamma_{\mcF^\heartsuit}$ may be identified with  $\mr{Hom}_{\mcO_{X^{(1)}}} (\mr{Ker}(\nabla_\mcF), \mr{Coker}(\gamma_{\mcF^\heartsuit}))$; this  is restricted to an identification between
the tangent space of $\mr{Quot}^\mr{triv} \left(\subseteq \mr{Quot}^0 \right)$ at the same point and $\mr{Ker}(\delta_{\mcF^\heartsuit})$.
Hence, $\mfO \mfp_{\mfs \mfl_n, g, 0}^{^\mr{Zzz...}} \times_{\mfM_{g,0}, x} k   \left(\cong \mr{Quot}^\mr{triv} \right)$ is unramified at the point classifying $\mcE^\spadesuit$ if and only if the equality $\mr{Ker}(\delta_{\mcF^\heartsuit}) =0$ holds.
This completes the proof of this assertion.
\end{proof}

\SSP
\bpr \label{wP109}
Suppose that  $X$ is smooth and of genus $0$.
Let us choose an $n$-theta characteristic $\vartheta := (\Theta, \nabla_\vartheta)$ of $X^\mr{log}/k$ (cf. Remark \ref{wR110}), and denote by $\mcF^\heartsuit := (\mcF, \nabla, \{ \mcF^j \}_j)$ the unique (up to isomorphism)  dormant $(\mr{GL}_n, \vartheta)$-oper with $\mcF^{\heartsuit \Rightarrow \spadesuit} \cong \mcE^\spadesuit$ (cf. (\ref{wE134})).
Then, $\mcE^\spadesuit$ is ordinary if the equality $H^1 (X^{(1)}, \mcE nd_{\mcO_{X^{(1)}}} (\mr{Ker}(\nabla))) =0$ holds.
\epr
\begin{proof}
The assertion follows from ~\cite[Chap.\,8, Propositions 8.4, 8.5, 8.11,  and 8.12]{Wak5}.
\end{proof}

Next, suppose that
$\msX$ is obtained by gluing together  pointed stable curves $\{ \msX_j \}_{j=1}^l$ ($l \geq 1$) by means of   clutching data (cf. ~\cite[\S\,7.1]{Wak5}).
Then,  for each $j=1, \cdots, l$,
 $\mcE^\spadesuit$ induces, via  restriction (and normalization in the sense of ~\cite[\S\,2.5, Proposition 2.19]{Wak5}),
 a dormant $\mfs \mfl_n$-oper $\mcE_j^\spadesuit$ on $\msX_j$ (cf. ~\cite[\S\,7.3, Proposition 7.12]{Wak5}).
 We shall refer to $\mcE^\spadesuit_j$ the {\bf restriction} of $\mcE^\spadesuit$ to $\msX_j$.
 
\SSP
\bpr \label{ghu2}
Let us keep the above notation.
Then, $\mcE^\spadesuit$ is ordinary  if and only if 
 $\mcE^\spadesuit_j$ is ordinary for every $j =1, \cdots, l$. 
 \epr
\begin{proof}
The assertion follows from  ~\cite[Chap.\,7, \S\,7.3, Theorem 7.13]{Wak5}.
\end{proof}

\LSP
\subsection{First proof of Theorem \ref{y019}, (i)} \label{S28}
We prove Theorem \ref{y019}, (i),  by using Proposition \ref{wP1231} described in the previous subsection.

With the notation in the statement of Theorem \ref{y019},
suppose further that the pull-back $w^*(\mcE^\spadesuit):= (w^*(\mcE), w^*(\nabla), w^*(\mcE_B))$ (where $\mcE^\spadesuit := (\mcE, \nabla, \mcE_B)$) is ordinary.
Let us consider the natural short exact sequence of complexes
\begin{align}
0 \migi \Omega_{X^\mr{log}/k^\mr{log}}\otimes (\mfs \mfl_n)_{\mcE_B}^{-1}  [-1] \migi \mcK^\bullet [\nabla^{\mr{ad}(0)}] \migi (\mfs \mfl_n)_{\mcE_B}^{0}[0] \migi 0,
\end{align}
where for an integer $m$ and a sheaf $\mcG$ we denote by $\mcG [m]$ the complex defined to be $\mcG$ concentrated at degree $-m$; this sequence induces 
 the sequence of $k$-vector spaces
\begin{align} \label{wE1598}
\Gamma (X, \Omega_{X^\mr{log}/k^\mr{log}} \otimes (\mfs \mfl_n)_{\mcE_B}^{-1}) \xrightarrow{q_{\mcE^\spadesuit}^\flat} \mbH^1 (X, \mcK^\bullet [\nabla^{\mr{ad}(0)}]) \xrightarrow{q_{\mcE^\spadesuit}^\sharp} H^1(X, (\mfs \mfl_n)_{\mcE_B}^0).
\end{align}
Note that
the pull-back of $((\mfs \mfl_n)_{\mcE_B}, \nabla^{\mr{ad}})$ by $w$ is  isomorphic to 
$((\mfs \mfl_n)_{w^*(\mcE_B)}, w^*(\nabla)^\mr{ad})$, and 
  the morphism 
$(\mfs \mfl_n)_{\mcE_B} \migi w_*((\mfs \mfl_n)_{w^*(\mcE_B)})$ induced by $w$ preserves the filtrations.
Hence,  we obtain  the following commutative diagram:
\begin{align} \label{fff11}
\vcenter{\xymatrix@C=36pt@R=36pt{
\Gamma (X, \Omega_{X^\mr{log}/k^\mr{log}} \otimes (\mfs \mfl_n)_{\mcE_B}^{-1})\ar[r]^-{q^\flat_{\mcE^\spadesuit}} \ar[d]&  \mbH^1 (X, \mcK^\bullet [\nabla^{\mr{ad}(0)}]) \ar[r]^-{q_{\mcE^\spadesuit}^\sharp} \ar[d]& H^1 (X, (\mfs \mfl_n)^{0}_{\mcE_B})  \ar[d] 
\\
\Gamma (Y, \Omega_{Y^\mr{log}/k^\mr{log}} \otimes (\mfs \mfl_n)_{w^*(\mcE_B)}^{-1}) \ar[r]_-{q_{w^*(\mcE^\spadesuit)}^\flat} & \mbH^1 (Y, \mcK^\bullet [w^*(\nabla)^{\mr{ad}(0)}])  \ar[r]_-{q_{w^*(\mcE^\spadesuit)}^\sharp} &  H^1 (Y, (\mfs \mfl_n)^{0}_{w^*(\mcE_B)}), 
}}
\end{align}
where all the vertical arrows arises from pull-back by $w$.
Here, recall (cf. ~\cite[Chap.\,2, \S\,2.1.4]{Wak5}) that, 
for each $j \in \mbZ$, 
 the subquotients $(\mfs \mfl_n^{j}/\mfs \mfl_n^{j+1})_{\mcE_B}$ and $(\mfs \mfl_n^{j}/\mfs \mfl_n^{j+1})_{w^*(\mcE_B)}$  decompose into  the direct sums of finite copies of $\Omega_{X^\mr{log}/k}^{\otimes j}$ and  $\Omega_{Y^\mr{log}/k}^{\otimes j}$ respectively.
Also,  the morphism $(\mfs \mfl_n^j/\mfs \mfl_n^{j+1})_{\mcE_B} \migi w_*((\mfs \mfl_n^j/\mfs \mfl_n^{j+1})_{w^*(\mcE_B)})$ induced by $w$ is compatible, via these decompositions,  with the natural morphism $\Omega_{X^\mr{log}/k} \migi w_*(\Omega_{Y^\mr{log}/k})$.
By using these facts, we see 
  that the leftmost and rightmost vertical arrows in (\ref{fff11}) are injective.
Hence, the middle arrow $\mbH^1 (X, \mcK^\bullet [\nabla^{\mr{ad}(0)}]) \migi \mbH^1 (Y, \mcK^\bullet [w^*(\nabla)^{\mr{ad}(0)}])$ turns out to be  injective.
This injection fits into the following commutative square diagram:
\begin{align} \label{fff11}
\vcenter{\xymatrix@C=46pt@R=36pt{
\mbH^1 (X, \mcK^\bullet [\nabla^{\mr{ad}(0)}])\ar[r]^-{\xi_{\mcE^\spadesuit}} \ar[d] & \Gamma (X, \mr{Coker}(\nabla^\mr{ad}))\ar[d]
\\
\mbH^1 (Y, \mcK^\bullet [w^*(\nabla)^{\mr{ad}(0)}]) \ar[r]_-{\xi_{w^*(\mcE^\spadesuit)}} & \Gamma (Y, \mr{Coker}(w^*(\nabla)^\mr{ad}))
}}
\end{align}
(cf. (\ref{wE33451}) for the defintion of $\xi_{(-)}$), where the right-hand vertical arrow arises naturally  from pull-back by $w$.
The lower horizontal arrow is injective because of the ordinariness assumption on $w^*(\mcE^\spadesuit)$ together with the result of Proposition \ref{wP1231}.
Hence, the commutativity of (\ref{fff11}) implies that $\xi_{\mcE^\spadesuit}$ is injective.
By Proposition  \ref{wP1231} again, $\mcE^\spadesuit$ turns out to be ordinary.
This completes the proof of Theorem \ref{y019}, (i).

\LSP
\subsection{Second proof of Theorem \ref{y019}, (i)} \label{Sf28}
In this subsection, we shall give another proof of Theorem \ref{y019}, (i), using Proposition \ref{wP3} under the assumption  that both  $\msX$ and $\msY$ are unpointed.
 
Let us choose a theta characteristic $\varTheta$ on $X$ and  a dormant $(\mr{GL}_n, \varTheta)$-oper $\mcF^\heartsuit := (\mcF, \nabla, \{ \mcF_j \}_j)$ on $X$ with
$\mcF^{\heartsuit \Rightarrow \spadesuit} \cong \mcE^\spadesuit$.
Denote by $w^{(1)}$ the base-change of $w$ by $F_{\mr{Spec}(k)}$; this morphism fits into the following commutative square diagram:
\begin{align} \label{DD0111}
\vcenter{\xymatrix@C=46pt@R=36pt{
Y \ar[r]^w \ar[d]_{F_{Y/k}} & X \ar[d]^{F_{X/k}} \\
Y^{(1)} \ar[r]_{w^{(1)}} & X^{(1)}.
}}
\end{align}
In particular,  this diagram induces an $\mcO_{Y^{(1)}}$-linear morphism
\begin{align} \label{E5002}
w^{(1)* } (F_{X/k*}(\mcF)) \migi F_{Y/k*} (w^* (\mcF)).
\end{align}
Here, let us prove the following lemma.

\SSP
\ble \label{P0d13d}
The morphism (\ref{E5002}) is an isomorphism, and restricted to
an isomorphism of $\mcO_{Y^{(1)}}$-modules
\begin{align} \label{EE01}
w^{(1)* } (\mr{Ker} (\nabla)) \isom \mr{Ker} (w^*(\nabla)).
\end{align}
\ele
\begin{proof}
The former assertion may be immediately verified  because
the   \'{e}taleness of $w$ implies that   (\ref{DD0111})
is  cartesian.

Next, let us consider the latter assertion.
Similarly to  (\ref{E5002}), we have the composite isomorphism
\begin{align} \label{E5001}
w^{(1)* } (F_{X/k*}(\Omega_{X/k} \otimes \mcF)) \isom & \ F_{Y/k*} (w^* (\Omega_{X/k}\otimes \mcF)) 
\isom F_{Y/k*} (\Omega_{Y/k} \otimes w^*(\mcF)),
\end{align}
where the second arrow arises from the \'{e}taleness of $w$.
The square diagram
\begin{align} \label{E5004}
\vcenter{\xymatrix@C=66pt@R=36pt{
w^{(1)* } (F_{X/k*}(\mcF)) \ar[r]^-{w^{(1)*}(F_{X/k*}(\nabla))} \ar[d]^-{\wr}_-{(\ref{E5002})}  & w^{(1)* } (F_{X/k*}(\Omega_{X/k} \otimes \mcF)) 
  \ar[d]_-{\wr}^-{(\ref{E5001})}
\\
F_{Y/k*} (w^* (\mcF))  \ar[r]_-{F_{Y/k*} (w^*(\nabla))} & F_{Y/k*} (\Omega_{Y/k}\otimes w^*(\mcF))
}}
\end{align}
is commutative.
In particular, 
by taking the kernels  of  the upper and  lower  horizontal arrows, we obtain   (\ref{EE01}), as desired.
\end{proof}
\SSP

Let us consider the composite isomorphism 
\begin{align} \label{wE120}
w^{(1)*}(F_{X/k*}(\mcF/\mcF^1)) \isom F_{Y/k*}(w^*(\mcF/\mcF^1))\isom   F_{Y/k*} (w^*(\mcF)/w^*(\mcF^1)),
\end{align}
where the first arrow is (\ref{E5002}) with $\mcF$ replaced by $\mcF/\mcF^1$.
This morphism makes the following square diagram commute:
\begin{align}
\vcenter{\xymatrix@C=46pt@R=36pt{
w^{(1)*}(\mr{Ker}(\nabla)) \ar[r]^-{w^{(1)*}(\gamma_{\mcF^\heartsuit})} \ar[d]_-{(\ref{EE01})}^-{\wr} & w^{(1)*}(F_{X/k*}(\mcF/\mcF^1)) \ar[d]_-{\wr}^{(\ref{wE120})}
\\
\mr{Ker} (w^*(\nabla))  \ar[r]_-{\gamma_{w^*(\mcF^\heartsuit)}} &  F_{Y/k*} (w^*(\mcF)/w^*(\mcF^1)).
}}
\end{align}
Hence,  this diagram induces   an isomorphism
\begin{align} \label{wE125}
w^{(1)*}(\mr{Coker}(\gamma_{\mcF^\heartsuit})) \isom \mr{Coker}(\gamma_{w^*(\mcF^\heartsuit)}).
\end{align}
Under the identifications given by (\ref{EE01}) and (\ref{wE125}), 
 the following square diagram is commutative:
\begin{align} \label{wE22}
\vcenter{\xymatrix@C=46pt@R=36pt{
\mr{Hom}(\mr{Ker}(\nabla), \mr{Coker}(\gamma_{\mcF^\heartsuit}))\ar[r]^-{\delta_{\mcF^\heartsuit}} \ar[d] & H^1 (X^{(1)}, \mcO_{X^{(1)}}) \ar[d] \\
\mr{Hom}(\mr{Ker}(w^*(\nabla)), \mr{Coker}(\gamma_{w^*(\mcF^\heartsuit)})) \ar[r]_-{\delta_{w^*(\mcF^\heartsuit)}} & H^1 (Y^{(1)}, \mcO_{Y^{(1)}})
}}
\end{align}
(cf. (\ref{wE23}) for the definition of $\delta_{(-)}$), 
where the both sides of vertical arrows are induced by pull-back by $w$.
The lower horizontal arrow $\delta_{w^*(\mcF^\heartsuit)}$ is injective because of 
 the ordinariness assumption on  $w^*(\mcE^\spadesuit)$ together with the result of Proposition \ref{wP3}.
 Hence the commutativity of  (\ref{wE22}) implies that 
   $\delta_{\mcF^\heartsuit}$  is   injective.
By Proposition \ref{wP3} again,  $\mcE^\spadesuit$ turns out to be  ordinary.
This completes the proof of Theorem \ref{y019}, (i).

\vspace{10mm}
\section{Cyclic log \'{e}tale  coverings of  pointed curves of low genus} \SSP

In this section, 
we consider dormant $\mfs \mfl_n$-opers on 
pointed smooth curves of genus $\leq 1$.
The goal of this section is to prove (cf. Proposition \ref{ghjk00098}) that the pull-back of any dormant $\mfs \mfl_n$-oper  on 
a generic genus-$1$ curve 
  by  a cyclic  log \'{e}tale covering (whose degree is prime to $p$)    is always ordinary.
This assertion will be a key ingredient in the proof of Theorem \ref{y019}, (ii).
Let $n$ be as in the previous section.

\LSP
\subsection{Vector bundles on  the  projective line} \label{S31}

Recall the Birkhoff-Grothendieck's theorem,  which asserts  that for any rank $n$ vector bundle $\mcG$ on the projective line  $\mbP$  over $k$   is isomorphic to a direct sum of $n$ line bundles
$\bigoplus_{j= 1}^n \mcO_{\mbP} (w_j)$,
where $w_{j_1} \leq w_{j_2}$ if $j_1 < j_2$.
The nondecreasing sequence  of integers $(w_j)_{j=1}^n$  depends only on the isomorphism class of the vector bundle $\mcG$.

\SSP
\bde \label{wD04}
 With the above notation, we shall say that $\mcG$ is {\bf of type $(w_j)_{j=1}^n$}.
Also, we shall say that $\mcG$ is {\bf of homogeneous type} if it is of type $(w_j)_{j=1}^n$ with $w_1 = w_n$ (i.e., $w_1 = w_2 = \cdots = w_n$).
 \ede
\SSP

\ble \label{6234}
 Let $s$ be an integer and  $\{ \mcG_l\}_{l\in \mbZ_{\geq 0}}$  a collection  of rank $n$ vector bundles on $\mbP$
 such that  $\mcG_l$ (for each $l$) is of degree $s+l$  and type $(w_{l, j})_{j=1}^n$ 
 (hence $\sum_{j=1}^n w_{l, j} = s+l$) with   $w_{l, n} - w_{l, 1} \leq 1$.
 Also, assume that we are given a  sequence of $\mcO_{\mbP}$-linear injections 
 \begin{align}
 \mcG_0 \migiincl \mcG_1 \migiincl \mcG_2 \migiincl \mcG_3 \migiincl \cdots. 
 \end{align}
  Then, there exists $l_0 \in \mbZ_{\geq 0}$ such that $\mcG_{l_0}$ is of homogeneous type. 
 \ele
\begin{proof}
The proof of this assertion is entirely similar to the proof of ~\cite[Chap.\,8, \S\, 8.4, Lemma 8.16]{Wak5}.
\end{proof}

\LSP
\subsection{A log \'{e}tale covering  of the $3$-pointed projective line} \label{S311}

For each value  $t  \in k \sqcup \{ \infty \}$, denote by $[t]$ the $k$-rational point of $\mbP$ determined by $t$.
Let $\LL$ be a positive  integer 
  prime to  $p$,  and  denote by  $\mu_\LL := \{   \zeta_1, \cdots, \zeta_\LL \}$ the group of $\LL$-power roots of unity in $k$.
 The collection of data
\begin{equation}
\msP^{(\LL+2)}_{} := (\mbP/k,  \{ [0], [\zeta_1], \cdots, [\zeta_\LL], [\infty] \})
\end{equation}
forms a $(\LL+2)$-pointed smooth genus-$0$ curve.
In particular, 
$\msP^{3}_{} \left(= (\mbP/k, \{ [0], [1], [\infty]\})\right)$
  is a unique (up to isomorphism) $3$-pointed stable curve of genus $0$.
We shall write 
\begin{align}
\mbP^{(\LL +2)\text{-}\mr{log}} := \mbP^{\msP^{(\LL+2)}\text{-}\mr{log}}.
\end{align}

\SSP
\bpr \label{500987} Let $\mcF^\heartsuit := (\mcF, \nabla, \{ \mcF^j \}_{j})$ be a dormant $\mr{GL}_n$-oper on $\msP^{(\LL+2)}$, and suppose that the vector bundle $\mr{Ker}(\nabla)$ on $\mbP^{(1)}$  is of type  $(w_j)_{j=1}^n$ (where $w_1 \leq \cdots \leq w_n$).
Then, we have the inequality
$w_n -w_1 \leq \LL +1$, or equivalently, $| w_i -w_j | \leq \LL +1$ for any $i$, $j \in \{1, \cdots, n \}$.
 \epr
\begin{proof}
The proof   is entirely similar to the proof of  ~\cite[Chap.\,8, \S\,8.4,  Proposition 8.17]{Wak5} after replacing  the $3$-pointed projective line ``$\msP$"  with   $\msP^{(\LL +2)}$.
\end{proof}
\SSP

Consider the endomorphism $\pi_\LL : \mbP \migi \mbP$ of $\mbP \left(= \mr{Proj}(k[x_1, x_2]) \right)$ corresponding to the $k$-algebra  endomorphism of $k[x_1, x_2]$ given by assigning $x_1 \mapsto x_1^\LL$ and $x_2 \mapsto x_2^\LL$.
This morphism 
is unramified over $\mbP \setminus \{ [0],[\infty]\}$ and satisfies that
$\pi^{-1}([0]) = \{[0]\}$, $\pi^{-1}([\infty]) = \{[\infty]\}$, and $\pi^{-1}([1]) = \{[\zeta_1], \cdots, [\zeta_\LL]\}$.
Also, it extends to a log \'{e}tale Galois covering 
\begin{equation} \label{DE01}
\pi_\LL^\mr{log} : \mbP^{(\LL +2)\text{-} \mr{log}} \migi \mbP^{3\text{-} \mr{log}}
\end{equation}
over $k$  with Galois group $\mu_\LL$.

\SSP
\ble\label{h01} Let $\mcF^\heartsuit := (\mcF, \nabla, \{ \mcF^j \}_j)$ be
a dormant $\mr{GL}_n$-oper on $\msP^{(3)}$.
Then, 
the natural morphism 
\begin{align} \label{ghhk}
\pi_\LL^{(1)*} (\mr{Ker} (\nabla)) \migi  \mr{Ker} (\pi_\LL^*(\nabla))
\end{align}
 is an isomorphism of $\mcO_{\mbP^{(1)}}$-modules.
 (This assertion remains true even when $(\mcF, \nabla)$ is replaced by  an arbitrary vector bundle equipped with a connection with vanishing $p$-curvature.)
 \ele
\begin{proof}
The morphism  $\pi^\mr{log}_\LL$ is \'{e}tale on  
 $\mbP\setminus \{ [0], [\infty] \}$, 
so it suffices to shows that (\ref{ghhk})  is an isomorphism on the formal neighborhoods of $[0]$ and $[\infty]$.
By similarity of  the argument,  we only consider the situation around $[0]$.
Let us fix 
 a local parameter  $t$ of $\mbP$ at $[0]$.
 Then, the formal neighborhood of $[0]$ in $\mbP$ (resp., $\mbP^{(1)}$) is isomorphic to $\mr{Spf}(k[\! [t]\! ])$ (resp., $\mr{Spf}(k[\! [t^p]\! ])$) and 
 the local parameter $t$ gives 
   an identification $\Gamma (\mr{Spf}(k[\! [t]\!]), \Omega_{\mbP/k} |_{\mr{Spf}(k[\! [t]\!])}) = k[\! [t]\! ]$.
It follows from ~\cite[Corollary 2.10]{O2}  
  that  the $t$-adic completion of $(\mcF, \nabla)$
is isomorphic to
 the direct sum  of line bundles with connection $(\mcO_{\mr{Spf}(k[\! [t]\! ])}, \nabla_a)$ for various $a \in \mbZ$ with $0 \leq a < p$, where $\nabla_a$ denotes a $k$-connection on  $\mcO_{\mr{Spf}(k[\! [t]\! ])}$ given by 
 $v \mapsto t \cdot \frac{\partial v}{\partial t} - a \cdot v$ for any $v \in k[\! [t]\! ]$.
 Hence,  the problem is reduced  to the case where $(\mcF, \nabla) = (\mcO_{\mr{Spf}(k[\! [t]\! ])}, \nabla_a)$ for such an integer $a$.
We shall identify, in a natural fashion,  the formal neighborhood of $[0]$ in the domain $\mbP$ of $\pi_\LL$ (resp., the domain $\mbP^{(1)}$ of $\pi_\LL^{(1)}$) with $\mr{Spf} (k[\! [t^{\frac{1}{\LL}}]\! ])$ (resp., $\mr{Spf} (k[\! [t^{\frac{p}{\LL}}]\! ])$).
The ring homomorphism   
$k[\! [t]\! ] \migi  k[\! [t^\frac{1}{\LL}]\! ]$ (resp., $k[\! [t^p]\! ] \migi  k[\! [t^\frac{p}{\LL}]\! ]$) obtained as the formal completion of  $\pi_\LL$ (resp., $\pi^{(1)}_\LL$) at $[0]$ is  given by $t \mapsto (t^{\frac{1}{\LL}})^\LL$ (resp., $t^p \mapsto (t^{\frac{p}{\LL}})^\LL$).
  Since  $\mr{Ker}(\nabla_a) |_{\mr{Spf}(k[\! [t^p]\! ])}$ is isomorphic to $\mcO_{\mr{Spf}(k[\! [t^p]\! ])} \cdot t^a$, we have
  \begin{align} \label{E602}
  \pi^{(1)*}_\LL (\mr{Ker}(\nabla_a)) |_{\mr{Spf}(k[[t^\frac{p}{\LL}]])}
  \cong 
  \mcO_{\mr{Spf}(k[\! [t^{\frac{p}{\LL}}]\! ])}
   \cdot t^a.
  \end{align}
   On the other hand, 
    $\pi_\LL^*(\nabla_a)  |_{\mr{Spf} (k[[t^\frac{1}{\LL}]])}$  is  given by $v \mapsto t^{\frac{1}{\LL}} \cdot \frac{\partial v}{\partial (t^{\frac{1}{\LL}})} - a \LL   \cdot v$ (for each $a \in k[\! [t^\frac{1}{\LL}]\! ]$).
  Hence, we have 
  \begin{align} \label{E601}
  \mr{Ker} (\pi^*_\LL (\nabla_a)) |_{\mr{Spf} (k[\! [t^\frac{p}{\LL}]\! ])} \cong 
  \mcO_{\mr{Spf}(k[\! [t^{\frac{p}{\LL}}]\! ])}
   \cdot (t^{\frac{1}{\LL}})^{a \LL}.
  \end{align}
  By (\ref{E602}) and (\ref{E601}), the morphism (\ref{ghhk}) turns out to be  an isomorphism on
  $\mr{Spf}(k[\! [t^{\frac{p}{\LL}}]\! ])$.
  This completes the proof of this lemma.
\end{proof}
\SSP

\ble \label{001}  For any  dormant $\mr{GL}_n$-oper $\mcF^\heartsuit$    on $\msP^{(3)}$,
 there exists a dormant $\mr{GL}_n$-oper $\mcG^\heartsuit : = (\mcG, \nabla, \{ \mcG^j \}_j)$ on $\msP^{(3)}$ equivalent to $\mcF^\heartsuit$ such that
the vector bundle  $\mr{Ker} (\nabla)$ on $\mbP^{(1)}$ is of homogeneous type.
 \ele
\begin{proof}
We shall   apply the discussion in the proof of  ~\cite[Chap.\,8, \S\,8.4, Proposition 8.18]{Wak5} to $\mcF^\heartsuit$.
(Note that $\mcF^\heartsuit$ defines a  dormant $(\mr{GL}_n, \vartheta)$-oper for some $n$-theta characteristic $\vartheta$.)
Then, we can obtain 
 a collection   $\{ \mcF^\heartsuit_m \}_{m \in \mbZ_{\geq 0}}$
 consisting of 
 dormant $\mr{GL}_n$-opers $\mcF^\heartsuit_m := (\mcF_m, \nabla_{m}, \{ \mcF_m^j \}_{j})$  on $\msP^{3}$ which  are all   equivalent to $\mcF^\heartsuit$ and satisfy $\mr{deg} (\mr{Ker}(\nabla_{m +1})) = \mr{deg} (\mr{Ker}(\nabla_{m})) + 1$.
 Moreover,  we can obtain  a sequence of $\mcO_{\mbP^{(1)}}$-linear injections
 \begin{align}
 \mr{Ker}(\nabla_0) \migiincl \mr{Ker}(\nabla_1) \migiincl \mr{Ker}(\nabla_2) \migiincl \cdots.
 \end{align}
 Now, suppose that there exists $m \geq 0$ for which  $\mr{Ker}(\nabla_{m})$ is of type $(w_{m, j})_{j=1}^n$ with $w_{m, n} -w_{m, 1} \geq 2$.
Let us take a positive integer $\LL$ with $\LL \geq 2$, $p \nmid \LL$  and take $\pi_\LL^\mr{log}$ as in (\ref{DE01}).
The vector bundle  $\mr{Ker} (\pi_\LL^*(\nabla_m))$ ($\cong \pi_\LL^{(1)*}(\mr{Ker}(\nabla_{m}))$   by Lemma \ref{h01})
 on $\mbP^{(1)}$ is of type $(\LL \cdot w_{m, j})_{j=1}^n$.
 By  Proposition \ref{500987}, we have
 \begin{equation}
 \left(\LL \cdot 2 \leq \right)   \LL\cdot w_{m, n} - \LL \cdot w_{m, 1} \leq \LL +1,
 \end{equation}
  which is a contradiction.
  This implies that, for any $m \geq 0$, $\mr{Ker}(\nabla_m)$
   is of type $(w_j)_{j =1}^n$ with $w_n -w_1 \leq 1$. 
Thus, 
by  Lemma  \ref{6234},
 there exists $m' \geq 0$ such that 
 $\mr{Ker}(\nabla_{m'})$ is of homogenous type.
 That is to say,  the dormant $\mr{GL}_n$-oper
  $\mcF^\heartsuit_{m'}$ satisfies the required conditions.
  \end{proof}
\SSP

By applying the above lemma, we obtain the following assertion.

\SSP
\bpr \label{PP01}  Let 
$\mcE^\spadesuit$
  be a dormant $\mfs \mfl_n$-oper  on $\msP^{(3)}$ (which is  ordinary, as mentioned in Remark \ref{wR135}).
 Then, the pull-back $\pi^*_\LL (\mcE^\spadesuit)$ of $\mcE^\spadesuit$  by  $\pi_\LL^\mr{log}$
 is ordinary.
 \epr
\begin{proof}
By the bijectivity of (\ref{wE134}) and Lemma \ref{001}, 
we can find a dormant $\mr{GL}_n$-oper  $\mcF^\heartsuit := (\mcF, \nabla,  \{ \mcF^j \}_j)$ on $\msP^{(3)}$ such that  
 $\mcF^{\heartsuit \Rightarrow \spadesuit} \cong \mcE^\spadesuit$ and $\mr{Ker}(\nabla)$  is of homogenous type.
It follows that the vector bundle $\mr{Ker} (\pi_\LL^*(\nabla))$ ($\cong \pi^{(1)*}_\LL (\mr{Ker} (\nabla))$ by Lemma \ref{h01})  is  of homogenous type, so
$\mcE nd_{\mcO_{\mbP^{(1)}}}  (\mr{Ker} (\pi_\LL^*(\nabla)))$ is isomorphic to a direct sum of finite copies of $\mcO_{\mbP^{(1)}}$.
This implies  
$H^1 (\mbP^{(1)}, \mcE nd_{\mcO_{\mbP^{(1)}}} (\mr{Ker} (\pi_\LL^*(\nabla)))) =0$.
Hence, the assertion follows from Proposition \ref{wP109}.
\end{proof}

\LSP
\subsection{A pointed totally degenerate curve} \label{S42}

Let us fix a positive integer $m$.
For each   positive integer $\LL$ prime to $p$,  we shall denote by 
\begin{equation}
\msC_{\LL} := (C_\LL/k , \{ \sigma_{\LL, i} \}_{i=1}^{ m \LL^2})
\end{equation}
the $m \LL^2$-pointed stable genus-$1$ curve over $k$  determined uniquely  by the following two conditions:
\begin{itemize}
\item
The underlying curve $C_\LL$ is given by 
$C_\LL = \bigcup_{j =1}^{m \LL } P_{\LL, j}$,
such that 
 $(P_{\LL, j}, \{ \sigma_{\LL, i} \}_{(j-1)\LL <  i \leq  j \LL}) = (\mbP, \{ [\zeta_1], \cdots, [\zeta_d]\})$ (where $\{ \zeta_1, \cdots, \zeta_d \} = \mu_d$) for every $j =1, \cdots m \LL$);
\item
There exist mutually distinct  $k$-rational points $\delta_{1}, \cdots, \delta_{m \LL}$  of  $C_m$
  such that 
\begin{equation}
P_{\LL, i} \cap P_{\LL, j} = \begin{cases}
\delta_{i} & \text{if $j = i+1$ mod $m \LL$} \\
\emptyset & \text{if otherwise.}
 \end{cases}
\end{equation}
(Hence, the set $\{ \delta_{1}, \cdots, \delta_{m \LL} \}$ coincides with the set of nodes of $C_\LL$).
\end{itemize}
In particular,
we obtain a log curve $C_\LL^\mr{log}/\mr{Spec}(k)^\mr{log}_\LL$, where $\mr{Spec}(k)^\mr{log}_\LL := \mr{Spec}(k)^{\msC_\LL \text{-}\mr{log}}$.
Moreover, 
 for each $j = 1, \dots, m \LL$,  we obtain (after ordering  the set of marked points) a $(\LL +2)$-pointed smooth genus-$0$ curve
\begin{equation}
\msP_{\LL, j} := (P_{\LL, j}, \{ \sigma_{\LL, i} \}_{(j-1)\LL < i\leq j \LL} \cup \{ \delta_{j-1}, \delta_{j} \})
\end{equation}
(where $\delta_{0} := \delta_{m \LL}$) over $\mr{Spec}(k)_{\LL, j} := \mr{Spec}(k)$.
It determines a log curve $P_{\LL, j}^\mr{log}/\mr{Spec}(k)_{\LL, j}^\mr{log}$.
Also, write $P_{\LL, j}^{\mr{log}^\dagger}$ for the log scheme obtained by equipping $P_{\LL, j}$ with the log structure pulled-back from $C_{\LL}^\mr{log}$ via the natural closed  immersion $P_{\LL, j} \migiincl C_{\LL}$.
The structure morphism $P_{\LL, j} \migi \mr{Spec}(k)$ extends to a morphism of log schemes  $P_{\LL, j}^{\mr{log}^\dagger} \migi \mr{Spec}(k)^\mr{log}_\LL$  (but it need not be a log curve).
Moreover, the natural morphisms $\mr{Spec}(k)^\mr{log}_{\LL} \migi \mr{Spec}(k)_{\LL, j}$ and $P_{\LL, j}^{\mr{log}^\dagger} \migi P_{\LL, j}$ extend to morphisms $\mr{Spec}(k)^\mr{log}_\LL \migi \mr{Spec}(k)_{\LL, j}^\mr{log}$ and $u^\mr{log}_{\LL, j} : P_{\LL, j}^{\mr{log}^\dagger} \migi P_{\LL, j}^\mr{log}$ respectively, which makes the following square diagram commute:
\begin{align}
\vcenter{\xymatrix@C=46pt@R=36pt{
P_{\LL, j}^{\mr{log}^\dagger} \ar[r]^-{u^\mr{log}_{\LL, j}} \ar[d] & P_{\LL, j}^\mr{log} \ar[d] \\
\mr{Spec}(k)^\mr{log}_\LL \ar[r] & \mr{Spec}(k)_{\LL, j}^\mr{log}.
}}
\end{align}
The morphism $\mcT_{P_{\LL, j}^{\mr{log}^\dagger}/\mr{Spec}(k)_\LL^\mr{log}} \isom u_{\LL, j}^*(\mcT_{P_{\LL, j}^\mr{log}/\mr{Spec}(k)_{\LL, j}^\mr{log}})$ induced by $u^\mr{log}_{\LL, j}$ is an isomorphism.
Hence, as mentioned in the comment preceding Proposition \ref{ghu2},
we can define (up to normalization) the restriction of each dormant $\mfs \mfl_n$-oper on $\msC_d$ to $P_{\LL, j}^\mr{log}/\mr{Spec}(k)_{\LL, j}^\mr{log}$.


\LSP
\subsection{Degeneration of a cyclic  covering of an elliptic curve} \label{S43}


Let   $v : \mr{Spec} (k) \migi \overline{\mfM}_{1, m}$ denote a $k$-rational point of the moduli stack $\overline{\mfM}_{1, m}$ classifying    $\msC_1$.
The divisor at infinity of $\overline{\mfM}_{1, m}$ around $v$ is locally given  as $\prod_{i=1}^m t_i =0$ for some local functions $t_i$ ($i = 1, \cdots, m$).
Then, the choices of such functions   determine
a morphism
$v_{R_1}  : \mr{Spec} (R_1) \migi \overline{\mfM}_{1, m}$, where $R_1:= k[[ t_1, \cdots, t_m]]$.
Denote by
\begin{align}
\msC_{1, R_1} := (C_{1, R_1}/R_1, \{ \sigma_{1, R_1, i} \}_{i =1}^{m})
\end{align}
the $m$-pointed  stable genus-$1$ curve  classified by $v_{R_1}$. 
In particular,  the special fiber is isomorphic to $\msC_1$ and the   underlying  curve of its  generic fiber is smooth.
The resulting log structure on  $\mr{Spec} (R_1)$ is the log structure   associated with the homomorphism
$\mbN^{\oplus m} \migi R_1$ given by assigning $(a_i)_{i=1}^m \mapsto \prod_{i=1}^m t_i^{a_i}$.

Next, let $\LL$ be a positive integer prime to $p$.
Write  $R_\LL :=  k [[ t^{\frac{1}{\LL}}_1, \cdots, t^{\frac{1}{\LL}}_m]]$  (resp., $R_\infty := \varinjlim_{p \nmid n} k [[ t^{\frac{1}{n}}_1, \cdots, t^{\frac{1}{n}}_m]]$) and equip  $\mr{Spec} (R_\LL)$ (resp., $\mr{Spec} (R_\infty)$) with the log structure associated with the homomorphism $(\frac{1}{\LL}\mbN)^{\oplus m} \migi R_\LL$  (resp., $\varinjlim_{p \nmid n} (\frac{1}{n}\mbN)^{\oplus m} \migi R_\infty$) given by $(a_i)_{i=1}^m \mapsto \prod_{i=1}^m t_i^{a_i}$.
Denote by  $\mr{Spec}(R_\LL)^\mr{log}$ (resp., $\mr{Spec}(R_\infty)^\mr{log}$)  the resulting log scheme and by 
$\mr{Spec} (k)^\mr{log}_\infty$   the log scheme defined as $\mr{Spec} (k)$ equipped with
the log structure pulled-back from  $\mr{Spec}(R_\infty)^\mr{log}$
via the closed immersion $\mr{Spec} (k) \migiincl  \mr{Spec}(R_\infty)$.
The sequence of homomorphisms  $R_1 \migiincl R_\LL \migiincl R_\infty \migisurj k$ yield a sequence of morphisms of log schemes
\begin{align} \label{EE013}
\mr{Spec} (k)^\mr{log}_\infty  \migi  \mr{Spec}(R_\infty)^\mr{log} \migi \mr{Spec}(R_\LL)^\mr{log} \migi \mr{Spec}(R_1)^\mr{log},
\end{align}
  by which $\mr{Spec}(R_\infty)^\mr{log}$ may be thought of as a universal log \'{e}tale covering of $\mr{Spec}(R_1)^\mr{log}$.

Let us fix an algebraically closed field $K$ together with an inclusion $R_\infty \migiincl K$, which induces a morphism $\mr{Spec}(K) \migi \mr{Spec} (R_\infty)^\mr{log}$.
Write 
\begin{align}
\msC_{1, K} := (C_{1, K}, \{ \sigma_{1, K,  i} \}_{i=1}^m)
\end{align}
for  the base-change of $\msC_{1, R_1}$ to $K$, 
and write
\begin{align}
C_{1, R_\infty}^{\mr{log}^\Box} := C_{1, R_1}^{\mr{log}^\Box}  \times_{\mr{Spec}(R_1)^{\mr{log}} } \mr{Spec}(R_\infty)^\mr{log},\end{align}
where $\Box$ denotes either the presence or absence of the prime ``$ \, '  \, $".
In particular, $\msC_{1, K}$ is classified by a geometric generic point of $\overline{\mfM}_{1, m}$ and 
we have $C_{1, K} \cong  C_{1, R_\infty}^{\mr{log}'} \times_{\mr{Spec}(R_\infty)^\mr{log}} \mr{Spec}(K)$.
The natural morphisms   $C_{1, K} \migi C_{1, R_\infty}^{\mr{log}'}$ and $C_{1, R_\infty}^{\mr{log}} \migi C_{1, R_\infty}^{\mr{log}'}$ give rise to a diagram  of logarithmic fundamental groups 
\begin{align} \label{DDD01}
\vcenter{\xymatrix@C=46pt@R=36pt{
&  \pi^{\mr{ket}}_1 (C_{1, R_\infty}^{\mr{log}})  \ar[d] & 
\\
\pi^{\mr{et}}_1 (C_{1, K})  \ar[r] &  \pi^{\mr{ket}}_1 (C_{1, R_\infty}^{\mr{log}'}) \ \Big( \ar[r]^{\hspace{-18mm}\sim} &  \pi^{\mr{ket}}_1 (C_1^{\mr{log}'} \times_{\mr{Spec}(k)_1^{\mr{log}}} \mr{Spec}(k)_\infty^\mr{log}) \Big),
}}
\end{align}
where all the arrows are determined up to   choices of base points, i.e., up to  composition with  inner automorphisms.
 The vertical arrow in (\ref{DDD01}) is surjective and
the horizontal arrow   becomes an isomorphism after taking maximal pro-$\LL$ quotient $(-)^{(\LL)}$.
Here,  let $[\LL ]_{K}$ denote the endomorphism of $C_{1,  K}$ determined by multiplication by $\LL$ (defined after fixing a base point of $C_{1,  K}$), which is an abelian \'{e}tale covering of degree $\LL^2$.
There exist a unique (up to isomorphism)  log scheme  $C_{\LL, R_\infty}^{\mr{log}}$  over $\mr{Spec}(R_\infty)^\mr{log}$ and   a log  \'{e}tale covering 
\begin{align}
[\LL]^\mr{log}_{R_\infty} : C_{\LL, R_\infty}^{\mr{log}} \migi C_{1, R_\infty}^{ \mr{log}}
\end{align}
 over $\mr{Spec} (R_\infty)^\mr{log}$  corresponding to   $[\LL]_{K}$
   via the resulting surjection
 \begin{align}
  \pi^{\mr{ket}}_1 (C_{1, R_\infty}^{\mr{log}})  \migisurj  \pi^{\mr{ket}}_1 (C_{1, R_\infty}^{\mr{log}'})^{(\LL)} \isom  \pi^{\mr{et}}_1 (C_{1, K})^{(\LL)}.
 \end{align} 
This means that $C_{\LL, R_\infty} \times_{R_\infty} K \cong C_{1,K}$ and  the fiber  of $[\LL]_{R_\infty}$ over $K$ coincides with $[\LL]_{K}$.
One may find 
an  $m \LL^2$-pointed stable genus-$1$ curve  
\begin{align}
\msC_{\LL, R_\LL} := (C_{\LL , R_\LL}/R_\LL, \{ \sigma_{\LL,  R_\LL, i} \}_{i=1}^{m \LL^2})
\end{align}
  over $R_\LL$ with $\mr{Spec}(R_\LL)^{\msC_{\LL, R_\LL} \text{-}\mr{log}} = \mr{Spec}(R_\LL)^\mr{log}$ (defined as above)  together with  a morphism $[\LL]^\mr{log}_{R_\LL} : C_{\LL, R_\LL}^{\mr{log}} \migi C_{1, R_1}^{\mr{log}}$
  which 
   makes  the following diagram  commute:
  \begin{align} \label{DDD01f}
\vcenter{\xymatrix@C=46pt@R=36pt{
C_{\LL, R_\infty}^{\mr{log}} \ar[r] \ar[d] &C_{\LL, R_\LL}^{\mr{log}} \ar[r]^{[\LL]^\mr{log}_{R_\LL}} \ar[d] & C_{1, R_1}^{\mr{log}}  \ar[d]
\\
\mr{Spec}(R_\infty)^\mr{log} \ar[r] & \mr{Spec}(R_\LL)^\mr{log}  \ar[r] & \mr{Spec}(R_1)^\mr{log},
}}
\end{align}
where the horizontal arrows are morphisms appearing in (\ref{EE013}), the left-hand square is  cartesian, and the base-change to $\mr{Spec}(R_\infty)^\mr{log}$ of the morphism $C_{\LL, R_\LL}^\mr{log} \migi C_{1, R_1}^\mr{log} \times_{\mr{Spec}(R_1)^\mr{log}} \mr{Spec}(R_\LL)^\mr{log}$ induced by  the right-hand square coincides with $[\LL]_{R_\infty}^\mr{log}$.
The special fiber of the right-hand square determines the following square diagram:
 \begin{align} \label{Dfg01f}
\vcenter{\xymatrix@C=46pt@R=36pt{
C_{\LL}^{\mr{log}} \ar[r]^{[\LL]_k^\mr{log}} \ar[d] & C_{1}^{\mr{log}}  \ar[d]
\\
 \mr{Spec}(k)_\LL^\mr{log}  \ar[r] & \mr{Spec}(k)_1^\mr{log}.
}}
\end{align}
The image  of each $P_{\LL, j} \subseteq C_{\LL}$ ($j = 1, \cdots, m \LL$) via $[\LL]_{k}$  is contained in
$P_{1, s_j} \subseteq C_{1}$ for some $s_j \in \{ 1, \cdots, m \}$.
Thus, the restriction of $[\LL]_{k}^\mr{log}$ to  $P_{\LL, j}$ defines a morphism  
\begin{align}
[\LL]_{k}^\mr{log} |_{P_{\LL, j}} : P_{\LL, j}^{\mr{log}^\dagger} 
 \migi P_{1, s_j}^{\mr{log}^\dagger}. 
\end{align}
By a straightforward argument, one may verify that  the square diagram
\begin{align} \label{wE4456}
\vcenter{\xymatrix@C=46pt@R=36pt{
P_{\LL, j}^{\mr{log}^\dagger} 
\ar[r]^-{[\LL]_{k}^\mr{log} |_{P_{\LL, j}}} \ar[d]_-{u_{\LL, j}^\mr{log}} & P_{1, s_j}^{\mr{log}^\dagger} 
 \ar[d]^-{u_{1, s_j}^\mr{log}}
\\
\mbP^{(\LL+2)\text{-}\mr{log}} 
 \ar[r]_-{\pi_\LL^\mr{log}} & \mbP_{}^{3\text{-}\mr{log}}
}}
\end{align}
(cf.  (\ref{DE01}) for the definition of $\pi^\mr{log}_\LL$)
is commutative.

\SSP
\ble \label{L001} 
Let $\mcE^\spadesuit$ be a dormant $\mfs \mfl_n$-oper on $C_1^\mr{log}/\mr{Spec}(k)_1^\mr{log}$ (i.e., on $\msC_1$).
Then, the pull-back $[\LL]_k^*(\mcE^\spadesuit)$ of $\mcE^\spadesuit$ by $[\LL]_k^\mr{log}$ is ordinary.
\ele
\begin{proof}
The assertion follows from  the commutativity of (\ref{wE4456}) and Propositions \ref{ghu2},  \ref{PP01}.
\end{proof}
\SSP

 The above  lemma implies  the following proposition.
 
\SSP
\bpr \label{ghjk00098} 
(Recall that $\msC_{1, K}$ is classified by  a geometric generic point of $\overline{\mfM}_{1, m}$.)
Let $(\msY, w)$ be a Galois \'{e}tale covering of $\msC_{1, K}$ over $K$ such that the Galois group $\mr{Gal}(\msY/\msC_{1, K})$ is abelian and satisfies $p \nmid \mr{Gal}(\msY/\msC_{1, K})$. 
Then, for any dormant $\mfs \mfl_n$-oper $\mcE^\spadesuit$ on $\msC_{1, K}$,
 its pull-back $w^*(\mcE^\spadesuit)$ by $w$ is ordinary.
 \epr
\begin{proof}
Write  $\msC_{\LL, K}$   for  the base-change of $\msC_{\LL, R_\LL}$  to $K$.
Note that $(\msY, w)$ is dominated by  a Galois \'{e}tale covering $(\msC_{\LL, K}, [\LL]_{K})$ for a sufficiently large $\LL$ (prime to $p$).
Hence,    by  Theorem \ref{y019}, (i),  it suffices to consider  the 
 case where   $(\msY, w) = (\msC_{\LL, K}, [\LL]_{K})$.
The fiber product  $\overline{\mfO} \mfp_{\mfs \mfl_n, 1,m}^{^\text{Zzz...}} \times_{\overline{\mfM}_{1, m}, v_{R_1}} \mr{Spec}(R_1)$ is isomorphic to a disjoint union of finite  copies of $\mr{Spec} (R_1)$.
It follows that
there exists a morphism $\widetilde{v}_{R_1} : \mr{Spec} (R_1) \migi \overline{\mfO} \mfp_{\mfs \mfl_n, 1,m}^{^\text{Zzz...}}$ whose composite with $\mr{Spec} (K) \migi \mr{Spec} (R_1)$  classifies  $\mcE^\spadesuit$.
Denote by 
 $\mcE^\spadesuit_{R_1}$ the dormant $\mfs \mfl_n$-oper on $\msC_{1, R_1}$ classified by $\widetilde{v}_{R_1}$ and by $\mcE_k^\spadesuit$ its special fiber.
 By Lemma \ref{L001}, the  special fiber of  the pull-back $[\LL]_{R_\LL}^*(\mcE^\spadesuit_{R_\LL})$, i.e., $[\LL]_{k}^*(\mcE_k^\spadesuit)$,  is ordinary.
 Therefore, since ${^\circledcirc}\overline{\mfO}\mfp_{\mfs \mfl_n, g,r}^{^\mr{Zzz...}}$ is open in $\overline{\mfO}\mfp_{\mfs \mfl_n, g,r}^{^\mr{Zzz...}}$,
  the pull-back $[\LL]_{R_\LL}^*(\mcE^\spadesuit_{R_\LL})$, as well as its base-change $[d]_K^*(\mcE^\spadesuit)$ to $K$,    is ordinary.
  This completes the proof of this proposition. 
 \end{proof}

\vspace{10mm}
\section{Pull-back of  ordinary dormant $\mfs \mfl_n$-opers by a cyclic  covering} \SSP


\LSP
\subsection{Proof of Theorem B, (ii)} \label{S44}

In this final section, we prove Theorem \ref{y019}, (ii), by applying the results obtained so far.
Let $\LLL$ be a positive integer with $p \nmid \LLL$.

\SSP
\bde \label{D001}
Let $\msX$ be a pointed stable curve over $k$.
We shall say that  $\msX$ is {\bf $(\LLL, n)$-dormant-ordinary}
 if 
 for any  dormant $\mfs \mfl_n$-oper  $\mcE^\spadesuit$ on $\msX$
 and any cyclic \'{e}tale covering  $(\msY, w)$ of $\msX$ whose degree divides $\LLL$,
  the pull-back $w^*(\mcE^\spadesuit)$ is  ordinary.
 \ede
\SSP

For each pair of nonnegative integers $(g, r)$ with $2g-2+r >0$, 
we shall denote by 
\begin{equation}
{^\circledcirc}\mfM_{g, r}^{\LLL, n}
\end{equation}
the locus of $\mfM_{g, r}$ classifying $(\LLL, n)$-dormant-ordinary curves.
Then, this locus forms an open substack of $\mfM_{g, r}$.
Indeed, for each positive integer  $\LLL'$ with $\LLL' \mid \LLL$, 
let us consider 
 the stack $\mfR_{g, r}^{\LLL'}$ classifying collections $(\msX, \msY, w)$ consisting of 
 an $r$-pointed smooth proper genus-$g$ curve $\msX$ over $k$ and a cyclic \'{e}tale covering $(\msY, w)$ of $\msX$ whose degree is $\LLL'$ (cf. ~\cite[\S\,1.2, Definition 1.1]{CF} for the case of $r=0$).
One may verify that $\mfR_{g, r}^{\LLL'}$
 is a Deligne-Mumford stack over $k$ and the projection $\kappa^{\LLL'} :  \mfR_{g, r}^{\LLL'} \migi \mfM_{g,r}$ given by $(\msX, \msY, w) \mapsto \msX$ is finite.
We shall 
write $\widetilde{\mfR}^{\LLL'} := \mfR_{g, r}^{\LLL'} \times_{\kappa_{\LLL'}, \mfM_{g,r}} \mfO \mfp_{\mfs \mfl_n, g, r}^{^\mr{Zzz...}}$ and write
\begin{align}
\xi^{\LLL'} : \widetilde{\mfR}^{\LLL'} \migi \mfO \mfp_{\mfs \mfl_n, (\LLL'-1)g +1, \LLL' r}^{^\mr{Zzz...}}
\end{align}
for the morphism  given by $(\msX, \msY, w), \mcE^\spadesuit) \mapsto (\msY, w^*(\mcE^\spadesuit))$.
Then, $\mfS^{\LLL'} := (\xi^{\LLL'})^{-1} (\mfO \mfp_{\mfs \mfl_n, (\LLL'-1)g +1, \LLL' r}^{^\mr{Zzz...}}\setminus {^\circledcirc}\overline{\mfO} \mfp_{\mfs \mfl_n, (\LLL'-1)g +1, \LLL' r}^{^\mr{Zzz...}})$
has a structure of (reduced) closed substack of $\widetilde{\mfR}^{\LLL'}$.
Since the natural projection $\widetilde{\kappa}^{\LLL'} : \widetilde{\mfR}_{\LLL'} \migi \mfM_{g,r}$ is 
finite, 
the  image  $\widetilde{\kappa}^{\LLL'} (\mfS^{\LLL'})$  of   $\mfS^{\LLL'}$ via this projection  forms a closed substack of $\mfM_{g,r}$.
But, it follows from the various definitions involved that
 the equality ${^\circledcirc}\mfM_{g,r}^{\LLL, n} = \mfM_{g, r} \setminus \bigcup_{\LLL' \mid \LLL}\widetilde{\kappa}^{\LLL'} (\mfS^{\LLL'})$ holds, so ${^\circledcirc}\mfM_{g, r}^{\LLL, n}$ forms an open substack of $\mfM_{g, r}$.
Hence, the proof of Theorem \ref{y019}, (ii), is now reduced to proving the following   assertion.

\SSP
\ble \label{L023} 
Any geometric generic point  $\eta : \mr{Spec}(L) \migi \mfM_{g, r}$ (where $L$ denotes  an algebraically closed field over $k$) of $\mfM_{g, r}$ lies in  ${^\circledcirc}\mfM_{g, r}^{\LLL, n}$.
In particular, ${^\circledcirc}\mfM_{g, r}^{\LLL, n}$ forms a dense open substack of $\mfM_{g,r}$.
 \ele
\begin{proof}
We only consider the former assertion because the latter assertion is a direct consequence of   the former assertion together with  the discussion preceding this lemma.
Moreover, since the case of $g=1$ has already been  proved in Proposition \ref{ghjk00098},
it suffices to consider the case of $g \geq 2$.
Denote by $\msC_L := (C_L/L, \{ \sigma_{L, i} \}_i)$
  the pointed  smooth curve  classified by $\eta$.
Let  us take a cyclic \'{e}tale covering $(\msC'_L, w_L)$ of $\msC_L$ of degree  $\LLL$ and a dormant $\mfs \mfl_n$-oper $\mcE^\spadesuit$ on $\msC_L$.
To conclude the assertion, we shall show that 
$w_L^*(\mcE^\spadesuit)$ is ordinary.

By applying the case of genus-$1$ curves, we can 
 take an $r$-pointed stable curve $\msC_0 := (C_0, \{ \sigma_{0, i}\}_i$  of genus $g$ over $k$ 
which is obtained by gluing together $(\LLL, n)$-dormant-ordinary smooth  genus-$1$ curves  $\{ \msD_{j} \}_{j=1}^g$ (where $\msD_j := (D_j, \{ \sigma_{j, i} \}_i)$)
and such that
 there exist mutually distinct $k$-rational points  $\delta^D_1, \cdots, \delta^D_{g-1}$ of $C_0$ satisfying 
\begin{equation}
D_{i} \cap D_{j} = \begin{cases}
\delta_{i}^D & \text{if $j = i+1$} \\
\emptyset & \text{if otherwise.}
 \end{cases}
\end{equation}
(Hence, the set $\{ \delta_{1}^D, \cdots, \delta_{g-1}^D \}$ coincides with the set of of nodes of $C_0$).
After possibly replacing $L$ with its extension field, 
we may assume that there exist an inclusion $R \migiincl L$ over $k$, where $R := k[[t_1, \cdots, t_{g-1} ]]$, and   an $r$-pointed  stable genus-$g$ curve $\msC_R := (C_R/R, \{ \sigma_{R, i}\}_i)$
over $R$ whose spacial fiber
  is isomorphic to 
$\msC_0$   and
    whose   fiber over $L$ is isomorphic to $\msC_L$.
It follows  from ~\cite[\S\,5, Lemma 5]{Nak} that there exists a  cyclic  \'{e}tale  covering 
$(\msC'_R, w_R)$ of $\msC_R$ over $R$
 whose fiber over $L$ is isomorphic to 
$(\msC_L, w_L)$.
Let us fix 
 $j \in \{1, \cdots, g\}$ and 
 choose  a connected  component $D'_j$ of $C'_0 \times_{C_0} D_j$, where $C'_0$ denotes the underlying curve of the spacial fiber of $\msC'_R$. 
Then,  the restriction of $(\msC'_R, w_R)$  to  $D'_j$ defines 
  a cyclic  \'{e}tale covering $(\msD'_j, w_{R, j})$ of $\msD_j$
  with $\mr{Gal}(\msD'_j/\msD_j) \mid \LLL$.
Observe that the finiteness of 
 the projection $\overline{\mfO} \mfp_{\mfs \mfl_n, g, r}^{^\mr{Zzz...}} \migi \overline{\mfM}_{g, r}$ implies that
$\mcE^\spadesuit$ extends to a dormant $\mfs \mfl_n$-oper $\mcE^\spadesuit_R$ on $\msC_R$.
Since $\msD_j$ is  $(\LLL, n)$-dormant-ordinary, 
the pull-back of $\mcE^\spadesuit$ to $\msD'_j$ is ordinary, so the spacial fiber
of $w_R^*(\mcE_R^\spadesuit)$ is ordnary (cf. Proposition \ref{ghu2}).
Hence, 
by the openness of ${^\circledcirc}\overline{\mfO} \mfp_{g,r}^{^\mr{Zzz...}}$,
the dormant $\mfs \mfl_n$-oper $w_R^*(\mcE^\spadesuit_R)$, as well as $w_L^*(\mcE^\spadesuit)$,  turns out to be ordinary.
This completes the proof of this assertion.
\end{proof}
\SSP

Consequently, we have proved Theorem \ref{y019}, (ii), which, in particular,  implies Theorem \ref{EW1}, (ii).

\SSP
\begin{rema}[Ordinariness for  $\mfg = \mfs \mfo_{2a +1}$, $\mfs \mfp_{2b}$] \label{wR21}
The result of Theorem \ref{y019}, (ii),  can be generalized to other Lie algebras $\mfg$, not just $\mfs \mfl_n$.
For example, suppose that $n$ is a positive integer with $p > 2n$ and  $\mfg$ is either $\mfs \mfo_{2a+1}$ with   $a = \frac{n-1}{2} \in \mbZ_{>0}$ or $\mfs \mfp_{2b}$ with  $b = \frac{n}{2} \in \mbZ_{>0}$.
In the same way as the case of $\mfs \mfl_n$, the ordinariness of dormant $\mfg$-opers is defined 
by means of the \'{e}tale (or equivalently, unramified) locus of the moduli stack $\overline{\mfO} \mfp_{\mfg, g,r}^{^\mr{Zzz...}}$ classifying pointed stable  curves equipped with a dormant $\mfg$-oper.
The natural inclusion $\mfg \migiincl  \mfs \mfl_n$ induces a closed immersion $\overline{\mfO} \mfp_{\mfg, g,r}^{^\mr{Zzz...}} \migiincl  \overline{\mfO} \mfp^{^\mr{Zzz...}}_{\mfs \mfl_n, g,r}$ (cf. ~\cite[Chap.\,5, \S\,5.4, Proposition 5.16]{Wak5}).
 It follows that any dormant $\mfg$-oper is ordinary whenever the associated dormant $\mfs \mfl_n$-oper is ordinary.
Hence,  by applying the generic \'{e}taleness of  $\overline{\mfO} \mfp_{\mfg, g,r}^{^\mr{Zzz...}}/\overline{\mfM}_{g,r}$ (cf. ~\cite[Theorem G]{Wak5})  and Theorem \ref{y019}, (ii), proved above, we see that the statement of Theorem \ref{y019}, (ii), remains true even when  $\mfs \mfl_n$ is replaced  by $\mfg$. 
\end{rema}



\vspace{10mm}


 
  \begin{thebibliography}{99}
\bibitem{ACGH}
D. Abramovich, Q. Chen, D. Gillam, Y. Huang, M. Olsson, M. Satriano and  S. Sun.
Logarithmic geometry and moduli.
\textit{Handbook of Moduli, Vol.\,I}, Adv. Lect. Math. {\bf 24}, Int. Press, Somerville, MA,  (2013), pp.\,1-61.



\bibitem{BK}
S. Bloch,  K. Kato,
$p$-adic \'{e}tale cohomology,
{\it  Inst. Hautes Etudes Sci. Publ. Math.} {\bf 63} (1986), pp. 107-152.


\bibitem{Bou}
I. I. Bouw,
The $p$-rank of curves and covers of curves.
{\it Courbes semi-stables et groupe fondamental en g\'{e}om\'{e}trie alg\'{e}brique,} Progr. Math. {\bf 187}  (2000), pp. 267-277.


\bibitem{CF}
A. Chiodo, G. Farkas, 
Singularities of the moduli space of level curves, 
{\it J. Eur. Math. Soc.}, {\bf 19} (2017), pp. 603-658.

\bibitem{DM}
P. Deligne and  D. Mumford,
The irreducibility of the space of curves of given genus.
{\it Publ. Math. I.H.E.S.} {\bf 36}  (1969), pp. 75-110.


\bibitem{Hos2}
Y. Hoshi,
Frobenius-projective structures on curves in positive characteristic.
{\it Publ. Res. Inst. Math. Sci. } {\bf 56} (2020), pp. 401-430.

\bibitem{IR}
L. Illusie, M.Raynaud,
Les suites spectrales associ\'{e}es au complexe de de Rham-Witt,
{\it  Inst. Hautes Etudes Sci. Publ.} {\bf 57} (1983), pp. 73-212.




\bibitem{JP}
K. Joshi and  C. Pauly,
Hitchin-Mochizuki morphism, opers and Frobenius-destabilized vector bundles over curves.
\textit{Adv. Math.} {\bf 274}, (2015), pp. 39-75.




\bibitem{KaFu}
F. Kato,
Log smooth deformation and moduli of log smooth curves.
\textit{Internat. J. Math.} {\bf 11}, (2000), pp. 215-232.

\bibitem{KaKa}
K. Kato,
Logarithmic structures of Fontaine-Illusie.
\textit{Algebraic analysis, geometry, and number theory}, John Hopkins Univ. Press, Baltimore, (1989), pp. 191-224.

 \bibitem{Kal}
N. M. Katz,
Nilpotent connections and the monodromy theorem: Applications of a result of Turrittin.
\textit{Inst. Hautes Etudes Sci. Publ. Math.} {\bf 39} (1970),  pp. 175-232.




\bibitem{Kn2}
F. F. Knudsen,
The projectivity of the moduli space of stable curves. II. The stacks $M_{g,r}$.
\textit{Math. Scand.} {\bf 52} (1983),  pp. 161-199.

\bibitem{Maz}
B. Mazur, 
Frobenius and the Hodge filtration, 
{\it Bull. Amer. Math. Soc.} {\bf 78} (1972), pp. 653-667.

\bibitem{Mzk1}
S. Mochizuki,
A theory of ordinary $p$-adic curves.
\textit{Publ. RIMS} {\bf 32} (1996),  pp. 957-1151.

\bibitem{Mzk2}
S. Mochizuki,
{\it Foundations of $p$-adic Teichm\"{u}ller theory}.
American Mathematical Society,  (1999).

\bibitem{Muk}
S. Mukai,
{\it  An Introduction to Invariants and Moduli}.
Cambridge studies in advanced mathematics, Cambridge University Press, (2003).



\bibitem{Nak}
S. Nakajima,
On generalized Hasse-Witt invariants of an algebraic curve.
\textit{Galois groups and their representations (Nagoya, 1981) (Y. Ihara, ed.), Adv. Stud. Pure Math.} {\bf 2}, North Holland/Kinokuniya, (1983),  pp. 69-88.



\bibitem{Og}
A. Ogus,
\textit{$F$-Crystals, Griffiths Transversality, and the Hodge Decomposition}.
Ast\'{e}risque {\bf 221}, Soc. Math. de France, (1994).



 \bibitem{O2}
B. Osserman,
Logarithmic connections with vanishing $p$-curvature.
\textit{J. Pure and Applied Algebra} {\bf 213} (2009),  pp. 1651-1664.


\bibitem{OP}
E. Ozman, R. Pries,
Ordinary and almost ordinary Prym varieties.
{\it Asian J. Math.} {\bf 23} (2019), pp. 455-478.

\bibitem{Ray1}
M. Raynaud,
Rev\^{e}tements des courbes en caract\'{e}ristique $p >0$ et ordinarit\'{e}.
\textit{Compositio Math.} {\bf 123} (2000),  pp. 73-88.


\bibitem{Wak}
Y. Wakabayashi,
An explicit formula for the generic number of dormant indigenous bundles.
\textit{Publ. Res. Inst. Math. Sci.} {\bf 50}  (2014),  pp. 383-409.


\bibitem{Wak2}
Y. Wakabayashi,
Spin networks, Ehrhart quasi-polynomials, and combinatorics of dormant indigenous bundles.
\textit{Kyoto J. Math.} {\bf 59} (2019), pp. 649-684.



\bibitem{Wak3}
Y. Wakabayashi,
The symplectic nature of the space of dormant indigenous bundles on algebraic curves.
\textit{arXiv: math. AG/1411.1197v3}, (2021).



\bibitem{Wak10}
Y. Wakabayashi,
Quantization on algebraic curves with Frobenius-projective structure,
\textit{arXiv: math. AG/2004.04283}, (2020).

\bibitem{Wak5}
Y. Wakabayashi,
A theory of dormant opers on pointed stable curves ---a proof of Joshi's conjecture---,
\textit{arXiv: math. AG/1411.1208v4}, (2021).

\bibitem{Wak11}
Y. Wakabayashi,
Frobenius-Ehresmann structures and Cartan geometries in positive characteristic,
\textit{arXiv: math. AG/1411.1208v4}, (2021).












\end{thebibliography}
\end{document}